# Stability conditions for the explicit integration of projection based nonlinear reduced-order and hyper reduced structural mechanics finite element models

C. Bach,[1,2] L. Song,[2] T. Erhart,[3] and F. Duddeck[1,4]

[1] *Technische Universität München, Munich, Germany*

[2] *BMW Group, Research and Innovation Centre, Munich, Germany*

[3] *Dynamore GmbH, Stuttgart, Germany*

[4] *Queen Mary, University of London, UK*

Projection-based nonlinear model order reduction methods can be used to reduce simulation times for the solution of many PDE-constrained problems. It has been observed in literature that such nonlinear reduced-order models (ROMs) based on Galerkin projection sometimes exhibit much larger stable time step sizes than their unreduced counterparts. This work provides a detailed theoretical analysis of this phenomenon for structural mechanics. We first show that many desirable system matrix properties are preserved by the Galerkin projection. Next, we prove that the eigenvalues of the linearized Galerkin reduced-order system separate the eigenvalues of the linearized original system. Assuming non-negative Rayleigh damping and a time integration using the popular central difference method, we further prove that the theoretical linear stability time step of the ROM is in fact *always* larger than or equal to the critical time step of its corresponding full-order model. We also give mathematical expressions for computing the stable time step size. Finally, we show that under certain conditions this increase in the stability time step even extends to some hyper-reduction methods. The findings can be used to compute numerical stability time step sizes for the integration of nonlinear ROMs in structural mechanics, and to speed up simulations by permitting the use of larger time steps.

## 1. INTRODUCTION

Projection-based nonlinear model order reduction methods can be used to reduce the simulation times of many problems described by partial differential equations (PDEs). While still an area of active research, such methods have already been applied to a wide range of problems. The speed-ups achievable by nonlinear model order reduction methods are known to be related to a number of different phenomena. Both implicit and explicit time integration methods benefit from hyper reduction procedures such as the Discrete Empirical Interpolation Method (DEIM), or the Energy Conserving Sampling and Weighting (ECSW) method. These methods are all associated with reducing the number of costly nonlinear function evaluations during the online simulations. The reduction of the dimension of the linear system of equations allows for a further speed-up in the case of implicit time integration schemes. This is not the case for explicit methods, where the integration is performed using only information already available from previous time steps.

However, it has been observed in literature that explicitly integrated reduced-order models can in fact benefit from considerable speed-ups due to larger stability time step sizes. This seems to have been initially reported by [1–3], and later by [4–6]. No detailed mathematical explanation or reasoning for this interesting phenomenon has been given in literature yet. The most closely related analysis for the case of modal truncation is given by [1]. This work offers a rigorous analysis of the numerical stability of explicitly integrated nonlinear reduced-order models equipped with a Galerkin projection. It is proven, for the first time as far as the authors are aware, that the stability time step of Galerkin ROMs with non-negative Rayleigh damping and a semi-orthogonal or $\mathbf{M}$-orthogonal [1] reduced basis, is in fact *always* larger than or equal to the equivalent stability time step of the corresponding full-order model for the central difference method, under conventional linear stability theory assumptions. That is, the reduced-order model inherits conditional stability as a property of the explicit time integration scheme, but the time step size necessary to stabilize the central difference time integration method is not smaller (and often larger) than the stability time step of the original problem. We show that this is fundamentally due to particular eigenvalue separation properties that arise naturally from the truncation of the reduced basis. We also give mathematical proofs that this phenomenon extends to some, but not all hyper reduced systems. For instance, DEIM is shown to be inherently less stable as it does not guarantee symmetry, positive definiteness or eigenvalue interlacing of the system matrices. By contrast, ECSW preserves symmetry and positive semi-definiteness, but not necessarily eigenvalue interlacing. Finally, we analyze the influence of the exact modal content of the reduced-order basis (ROB) $\mathbf{V}$ on the system eigenvalues and the stability time step.

The implication is that reduced-order models are naturally more stable numerically than their full-dimensional counterparts. This means that larger – or at least equal size – time steps can be used in the online phase of the solution of a dynamical reduced or hyper reduced model. We nevertheless also emphasize that overly large time steps for highly nonlinear problems may violate linear stability theory assumptions. Numerical stability is not to be confused with consistency. While reduced-order modelling can be beneficial for numerical stability of explicit integration methods, it typically introduces consistency errors that may negatively affect convergence.

---

[1] where $\mathbf{M}$ is a diagonal, positive matrix

**Notation.**
Matrices or matrix valued functions are denoted using bold capital letters. Vectors or vector valued functions are written in lower-case letters, also in bold. The spectral norm of a matrix $\mathbf{A}$ is denoted by $\|\mathbf{A}\|_2$, its Frobenius norm is written as $\|\mathbf{A}\|_F$. Time derivatives are denoted using Newton's notation, where $\frac{d}{dt}\mathbf{x} = \dot{\mathbf{x}}$. In line with most mathematical literature and to enable a better understanding, we denote the (real) eigenvalues $\lambda_i$ and singular values $\sigma_i$ of a matrix arranged in increasing order throughout this work. Notice that this is contrary to the decreasing order convention which is widely used in data analysis and nonlinear model order reduction.

## 2. NUMERICAL INTEGRATION OF THE EQUATIONS OF MOTION

The semi-discrete form of the nonlinear equations of motion describing a general structural dynamics problem can be written as

$$\mathbf{M}(\dot{\mathbf{x}}, \mathbf{x}, t)\ \ddot{\mathbf{x}} + \mathbf{C}(\dot{\mathbf{x}}, \mathbf{x}, t)\ \dot{\mathbf{x}} + \mathbf{K}(\dot{\mathbf{x}}, \mathbf{x}, t)\ \mathbf{x} = \mathbf{f}_{\text{ext}}(\dot{\mathbf{x}}, \mathbf{x}, t), \tag{1}$$

where $\mathbf{x} \in \mathbb{R}^m$ is the vector of nodal displacements of the spatial discretization, $\mathbf{M} \in \mathbb{R}^{m\times m}$ denotes the mass matrix, $\mathbf{C} \in \mathbb{R}^{m\times m}$ is the damping matrix, and $\mathbf{K} \in \mathbb{R}^{m\times m}$ is the nonlinear stiffness matrix. $\mathbf{f}_{\text{ext}} \in \mathbb{R}^m$ denotes the externally applied loads. The time integration of these equations can be performed either using implicit or using explicit methods.

Implicit integration methods solve the system of equations (1), ensuring equilibrium at discrete time instants using an iterative method (e.g. the Newton-Raphson algorithm). Implicit methods are unconditionally stable and do not impose any *stability* limits on the time step size, under linear stability theory assumptions. However, the convergence of this iterative procedure at every time step can be an issue when dealing with highly nonlinear phenomena such as contact. The accuracy of the linearization becomes poor in the presence of large time steps or highly nonlinear terms, thus violating linear stability theory assumptions. It is therefore generally very difficult to solve highly nonlinear problems using large time steps. Moreover, implicit methods require the inversion of the linearized stiffness matrix $\mathbf{K}$ to solve for $\mathbf{x}$. This is a costly operation, but for many cases the effort for the inversion is offset by the possibility to use large time steps. Examples include the Newmark beta method [7], or the Houbolt and Bathe methods [8].

Explicit methods do not actually solve a system of equations. Using a diagonal lumped mass matrix, the updated nodal displacements and velocities at time $t_n$ only depend on previously computed quantities [9]. The updates can therefore be computed directly, avoiding costly convergence iterations and stiffness matrix inversions. In fact, the global damping or stiffness matrices are usually never assembled as a whole, saving memory. Instead, the term $\mathbf{C}(\dot{\mathbf{x}}, \mathbf{x}, t)\ \dot{\mathbf{x}} + \mathbf{K}(\dot{\mathbf{x}}, \mathbf{x}, t)\ \mathbf{x}$ is often replaced by a vector of nonlinear *internal forces* $\mathbf{f}_{\text{int}}(\dot{\mathbf{x}}, \mathbf{x}, t)$, comprising nonlinear contributions from contact, plasticity or external loads. However, explicit integration schemes, such as the central difference method, are only conditionally stable and typically require small time steps. Explicit methods are often used for the simulation of highly nonlinear phenomena with large deformations, such as impact and crash analysis, or natural disaster simulation. The focus of this work is the numerical stability of explicit methods.

### 2.1. Central difference method

The central difference method is still the most popular explicit time integration scheme for structural mechanics problems. It is second-order accurate and both simple to understand and to implement efficiently in practice. It is part of many widely used finite element solvers such as LS-DYNA [10] and Abaqus/explicit [11]. This is why it serves as an example for the stability analysis of reduced-order systems throughout the remainder of this paper.

Algorithm 1 provides an overview of a typical numerical implementation for the integration of the equations of motion (1). The following analysis is furthermore restricted to Lagrangian meshes. For the sake of simplicity (but without loss of generality), we further assume a constant time step size $\Delta t = t^{n+1} - t^n$, where $t$ is the time, and the superscript $(\cdot)^n$ is used to denote quantities at the $n$-th discretization point in time. In practice, the time step size $\Delta t$ may of course be adjusted during a simulation (cf. [9]). Recalling equation (1), the central difference method computes the velocities $\dot{\mathbf{x}}$ at time $t^{n+1/2}$ as

$$\dot{\mathbf{x}}^{n+1/2} = \frac{\mathbf{x}^{n+1} - \mathbf{x}^{n}}{\Delta t}, \qquad \text{with} \qquad t^{n+1/2} = t^{n} + \frac{\Delta t}{2}. \tag{2}$$

The corresponding formula for the acceleration terms $\ddot{\mathbf{x}}$ at time $t^n$ is given by

$$\ddot{\mathbf{x}}^{n} = \frac{\dot{\mathbf{x}}^{n+1/2} - \dot{\mathbf{x}}^{n-1/2}}{\Delta t} = \frac{\mathbf{x}^{n+1} - 2\mathbf{x}^{n} + \mathbf{x}^{n-1}}{\Delta t^2}. \tag{3}$$

Notice that there is a half time step lag between the velocities and the displacements and accelerations. While this version of the central difference scheme is very popular [9, 10], slight modifications thereof are also possible (cf. [5, 8]). The accelerations are related to the forces by the mass matrix $\mathbf{M}$ such that

$$\mathbf{M}^{n}(\dot{\mathbf{x}}^{n-1/2}, \mathbf{x}^{n}, t^{n})\ddot{\mathbf{x}}^{n} = \mathbf{f}^{n}_{\text{node}}(\dot{\mathbf{x}}^{n-1/2}, \mathbf{x}^{n}, t^{n}), \tag{4}$$

where $\mathbf{f}_{\text{node}}$ is the vector of nodal forces acting on every degree of freedom (DoF) of the assembled system. $\mathbf{f}_{\text{node}}$ comprises contributions from internal and external forces. Recalling equation (1), it can be written as

$$\mathbf{f}^{n}_{\text{node}}(\dot{\mathbf{x}}^{n-1/2}, \mathbf{x}^{n}, t^{n}) = \mathbf{f}^{n}_{\text{ext}}(\dot{\mathbf{x}}^{n-1/2}, \mathbf{x}^{n}, t^{n}) - \mathbf{C}^{n}(\dot{\mathbf{x}}^{n-1/2}, \mathbf{x}^{n}, t^{n})\dot{\mathbf{x}}^{n-1/2} - \mathbf{K}^{n}(\dot{\mathbf{x}}^{n-1/2}, \mathbf{x}^{n}, t^{n})\mathbf{x}^{n}. \tag{5}$$

The procedure for the numerical integration of (1) using the central difference scheme is summarized in algorithm 1.

**Algorithm 1** Central difference integration scheme $(\cdot)^n \rightarrow (\cdot)^{n+1}$, adapted from [9, 12]

**Input:** Previous displacements $\mathbf{x}^n$ and velocities $\dot{\mathbf{x}}^{n-1/2}$, nodal forces $\mathbf{f}^{n}_{\text{node}}$, mass matrix $\mathbf{M}^n$, time $t^n$, time step size $\Delta t$. If $n = 0$, then use the initial conditions $\dot{\mathbf{x}}^0$ instead of $\dot{\mathbf{x}}^{n-1/2}$.
**Output:** Updated displacements $\mathbf{x}^{n+1}$, velocities $\mathbf{x}^{n+1/2}$, time $t^{n+1}$.

1: Accelerations update: $\ddot{\mathbf{x}}^{n} \leftarrow (\mathbf{M}^{n})^{-1}\mathbf{f}^{n}_{\text{node}}$
2: Velocity update: $\dot{\mathbf{x}}^{n+1/2} \leftarrow \dot{\mathbf{x}}^{n-1/2} + \Delta t\ddot{\mathbf{x}}^{n}$
3: Displacements update: $\mathbf{x}^{n+1} \leftarrow \mathbf{x}^{n} + \Delta t\dot{\mathbf{x}}^{n+1/2}$
4: Time update: $t^{n} \leftarrow t^{n} + \Delta t;\quad n \leftarrow n + 1$

Due to the velocity time lag of $\Delta t/2$, the discrete form of the equations of motion discretized by the central difference method is [9, p. 400]

$$\mathbf{M}^{n}\ddot{\mathbf{x}}^{n} + \mathbf{C}^{n}\dot{\mathbf{x}}^{n-1/2} + \mathbf{K}^{n}\mathbf{x}^{n} = \mathbf{f}^{n}_{\text{ext}}. \tag{6}$$

For the purpose of this study, we assume that it is known how the nodal forces $\mathbf{f}_{\text{node}}$ are computed from the elemental contributions, contact forces, and externally applied loads. Any matrices that depend on $(\dot{\mathbf{x}}, \mathbf{x}, t)$ are typically linearized at specific points in time (cf. section 3.2, though not necessarily at every time step. The mass matrix $\mathbf{M}$ can be assumed as constant for some applications, if it is not configuration dependent and no artificial dynamic mass scaling is used.

## 3. LINEAR STABILITY THEORY

There are various ways of defining stability, and the stability analysis of nonlinear dynamical systems is in general a challenging problem. Throughout this work, we employ a notion of stability as described in [9], which originates from the works of Lyapunov [13].

**Definition.** *Consider the semi-discrete form of a dynamical system described by some arbitrary evolution equation. Let* $\mathbf{x}^0$ *be a vector of initial conditions, and let* $\delta$ *be a perturbation vector, with* $0 < \|\delta\|_2 \leq \varepsilon$*. Let further* $\mathbf{x}_A(t)$ *denote the solution of the system with the initial conditions* $\mathbf{x}_A(t = 0) = \mathbf{x}^0$*, and let* $\mathbf{x}_B(t)$ *denote the solution of the perturbed system, with* $\mathbf{x}_B(t = 0) = \mathbf{x}^0 + \delta$*. The solution is then called stable if*

$$\|\mathbf{x}_A(t) - \mathbf{x}_B(t)\|_2 \leq C\varepsilon \qquad \forall t > 0 \tag{7}$$

*for some* $C > 0$*. A more detailed explanation of this notion of stability can be found in [9].*

Physically, the above definition means that points that are close to each other stay close to each other as time passes. Notice that some physical systems exhibit instabilities that violate this principle, such as bifurcation points. The following sections provide an overview of how the numerical stability of time integration methods can be assessed.

### 3.1. Stability of time integration methods for linear systems

Provided that a physical system is stable, one can derive conditions under which different time integration schemes are also stable. In line with definition 3, a numerical solution is called stable, if small perturbations $\delta$ of the initial conditions do not lead to divergent behaviour of the solutions as time passes, i.e. for all time instants $n$ and initial conditions $\mathbf{x}^0 = \mathbf{x}_A^0$ (cf. [9])

$$\|\mathbf{x}_A^n - \mathbf{x}_B^n\|_2 \leq C\varepsilon \qquad \text{such that } \|\delta\|_2 \leq \varepsilon. \tag{8}$$

The stability of time integrators is usually assessed for linear systems, using Dahlquist's test equation [14]

$$\dot{y}(t) = \alpha y(t) \qquad \text{with} \qquad \alpha \in \mathbb{C}, \tag{9}$$

which has the trivial solution $y(t) = \mathrm{e}^{\alpha t}$. Depending on the time integration scheme, different *stability regions* arise. For the simple example of the explicit forward Euler method, $\dot{y}(t)$ is approximated by

$$\dot{y}^n = \frac{y^{n+1} - y^n}{\Delta t}. \tag{10}$$

Inserting (10) into (9) yields

$$y^{n+1} = (1 + \alpha \Delta t) y^n. \tag{11}$$

It can be shown that the above scheme remains stable as long as $|1 + \alpha\Delta t| \leq 1$, which is only possible for $Re(\alpha) < 1$. The product of the characteristic value $\alpha$ of the equation and the time step $\Delta t$ is often denoted as $z$. The stability region of the explicit Euler method is a circle in the complex $z$ plane, i.e. there is a constraint on how large the time step $\Delta t$ can be selected in order to ensure stability. This observation leads to the two following definitions which are important for the remainder of this paper.

**Definition.** *(Conditional stability)*
*A time integration scheme that ensures numerical stability of the solution of a (stable) linear system only when the time step size $\Delta t$ fulfils a stability condition $\Delta t \leq \Delta t_{crit}$ is called conditionally stable. The problem dependent upper bound on the admissible time step size $\Delta t_{crit}$ is called the critical time step.*

**Definition.** *(Unconditional stability)*
*A time integration scheme is called unconditionally stable if it ensures numerical stability of the solution of a (stable) linear system regardless of the chosen time step size $\Delta t$.*

Notice that numerical stability is only defined for stable systems. If the physical system of interest does not fulfil the stability condition 3, the stability of numerical methods applied to such a problem cannot be assessed. Even if the integration method itself is considered as stable, condition (8) will still be violated [9]. However, the analysis of highly nonlinear and dynamical phenomena such as impact, automotive crash, or natural disasters is of great practical interest. The popular assumption (or hope) in these cases is that numerically stable methods will at least "behave well" when applied to unstable systems [9, p. 392].
The following sections briefly detail how the stability condition for the explicit central difference method in the context of explicit FEM is derived, and how the critical time step size can be estimated at low computational costs.

### 3.2. Linearization of nonlinear systems

The above definitions of conditional and unconditional stability, as well as the analysis of the stability region of the explicit Euler method, were restricted to linear systems. The nonlinear equations of motion (1) are therefore linearized around certain points in time before linear stability theory can be applied. This need not necessarily be done at every time step. However, highly nonlinear problems – e.g. involving contact – require more frequent linearizations. As long as the periods between linearizations are small enough, the application of linear stability methods to nonlinear problems is usually deemed acceptable.
The need to repeatedly linearize the equations during the integration has important consequences:

- otherwise unconditionally stable methods may not use arbitrarily large time steps, because they could otherwise become unstable and yield incorrect results. This is the case, for example, when unconditionally stable implicit time integration schemes are used in the presence of contact nonlinearities, e.g. impact simulation.

- the critical time step size for conditionally stable methods may be somewhat smaller than the one predicted under linear stability theory assumptions. A safety factor is often used in practice to account for such phenomena.

In the following we denote the linearized mass, damping and stiffness matrices as $\mathbf{M}, \mathbf{C}, \mathbf{K}$. The semi-discrete form of the linearized equations around an arbitrary point then read

$$\mathbf{M}\ddot{\mathbf{x}} + \mathbf{C}\dot{\mathbf{x}} + \mathbf{K}\mathbf{x} = \mathbf{f}_{\text{ext}}(\dot{\mathbf{x}}, \mathbf{x}, t). \tag{12}$$

### 3.3. Amplification matrix notation and modal equations

For linear (or linearized) problems, the numerical stability of integration methods is often evaluated by analyzing the eigenvalues of the so-called amplification matrix. Consider the unreduced equations of motion (1), discretized in time by a central difference scheme (6). Inserting (2) and (3) into the linearized equation of motion for $t = t^n$ yields

$$\mathbf{M}\left(\frac{\mathbf{x}^{n+1} - 2\mathbf{x}^n + \mathbf{x}^{n-1}}{\Delta t^2}\right) + \mathbf{C}\left(\frac{\mathbf{x}^n - \mathbf{x}^{n-1}}{\Delta t}\right) + \mathbf{K}\mathbf{x}^n = \mathbf{f}^n_{\text{ext}} \tag{13}$$

$$\iff \frac{\mathbf{M}}{\Delta t^2}\mathbf{x}^{n+1} + \left(\mathbf{K} - \frac{2\mathbf{M}}{\Delta t^2} + \frac{\mathbf{C}}{\Delta t}\right)\mathbf{x}^n + \left(\frac{\mathbf{M}}{\Delta t^2} - \frac{\mathbf{C}}{\Delta t}\right)\mathbf{x}^{n-1} = \mathbf{f}^n_{\text{ext}} \tag{14}$$

$$\iff \mathbf{x}^{n+1} = \Delta t^2 \mathbf{M}^{-1}\left[\mathbf{f}^n_{\text{ext}} - \left(\mathbf{K} - \frac{2\mathbf{M}}{\Delta t^2} + \frac{\mathbf{C}}{\Delta t}\right)\mathbf{x}^n - \left(\frac{\mathbf{M}}{\Delta t^2} - \frac{\mathbf{C}}{\Delta t}\right)\mathbf{x}^{n-1}\right]. \tag{15}$$

Let $\mathbf{y}$ denote an auxiliary state variable such that

$$\mathbf{y}^n = \begin{bmatrix}\mathbf{x}^n \\ \mathbf{x}^{n-1}\end{bmatrix}. \tag{16}$$

Then we can express $\mathbf{y}^{n+1}$ as

$$\mathbf{y}^{n+1} = \begin{bmatrix}\mathbf{x}^{n+1} \\ \mathbf{x}^n\end{bmatrix} = \underbrace{\begin{bmatrix}2\mathbf{I} - \Delta t^2\mathbf{M}^{-1}\mathbf{K} - \Delta t\mathbf{M}^{-1}\mathbf{C} & \Delta t\mathbf{M}^{-1}\mathbf{C} - \mathbf{I} \\ \mathbf{I} & \mathbf{0}\end{bmatrix}}_{\mathbf{A}} \begin{bmatrix}\mathbf{x}^n \\ \mathbf{x}^{n-1}\end{bmatrix} + \begin{bmatrix}\Delta t^2\mathbf{M}^{-1}\mathbf{f}^n_{\text{ext}} \\ \mathbf{0}\end{bmatrix} \tag{17}$$

$$= \mathbf{A}\mathbf{y}^n + \begin{bmatrix}\Delta t^2\mathbf{M}^{-1}\mathbf{f}^n_{\text{ext}} \\ \mathbf{0}\end{bmatrix}. \tag{18}$$

The matrix $\mathbf{A} \in \mathbb{R}^{2m\times 2m}$ is called the amplification matrix. Any errors in $\mathbf{y}^n$ remain bounded for $n \to \infty$ if

$$\|\mathbf{A}\|_2 = \rho(\mathbf{A}) = \max_i |\lambda^A_i| \leq 1, \tag{19}$$

under the additional condition that any eigenvalue $\lambda^A_i$ of $\mathbf{A}$ which has an absolute magnitude $|\lambda_i| = 1$ may not have an algebraic multiplicity larger than 1 [8]. As $\mathbf{A}$ is generally not symmetric, its eigenvalues $\lambda^A_i$ are complex numbers and the stability region can be visualized as a circle with unit radius, centered at the origin of the complex plane.

It is computationally very expensive (nevertheless possible) to evaluate the eigenvalues of the above amplification matrix $\mathbf{A}$ numerically. There are more efficient ways of estimating or bounding its eigenvalues by exploiting the block structure. A very common alternative, however, is to use the modal equations to analyze stability. Let $\mu_i, \mathbf{v}_i, 1 \leq i \leq m$ denote the $m$ eigenpairs of the matrix $\mathbf{M}^{-1}\mathbf{K}$. Since the $\mathbf{v}_i$ are linearly independent, it is possible to express the displacements and their derivatives using a linear combination thereof. This corresponds to a basis change [2].

$$\mathbf{x}(t) = \sum_{i=1}^{m} \tilde{x}_i(t)\mathbf{v}_i; \quad \dot{\mathbf{x}}(t) = \sum_{i=1}^{m} \dot{\tilde{x}}_i(t)\mathbf{v}_i; \quad \ddot{\mathbf{x}}(t) = \sum_{i=1}^{m} \ddot{\tilde{x}}_i(t)\mathbf{v}_i \tag{20}$$

[2] Notice the similarity to projection based MOR methods, where $\mathbf{x}$ is only *approximated* using a truncated sum from 1 to $k < m$ (cf. section 5).

Assuming Rayleigh damping, i.e. $\mathbf{C} = a_1\mathbf{M} + a_2\mathbf{K}$, the modal form of the equations (1) can be written as [9, p. 400]

$$\ddot{\tilde{x}} + (a_1 + a_2\mu_j)\dot{\tilde{x}} + \mu_j\tilde{x} = \mathbf{v}_j^T\mathbf{f}_{\text{ext}}, \qquad j \in \{1, 2, \ldots, m\}. \tag{21}$$

$$\iff \ddot{\tilde{x}} + 2\xi_j\sqrt{\mu_j}\dot{\tilde{x}} + \mu_j\tilde{x} = \mathbf{v}_j^T\mathbf{f}_{\text{ext}} \quad \text{with} \quad \xi_j = \frac{a_1}{2\sqrt{\mu_j}} + \frac{a_2\sqrt{\mu_j}}{2} \tag{22}$$

The above is obtained by exploiting the orthogonality properties of $\mathbf{v}_i$ with respect to $\mathbf{M}$ and $\mathbf{K}$ [9, p. 393]

$$\mathbf{v}_j^T\mathbf{M}\mathbf{v}_i = \delta_{i,j}, \qquad \mathbf{v}_j^T\mathbf{K}\mathbf{v}_i = \mu_i^2\delta_{i,j} \quad \text{(no summation)}, \tag{23}$$

where $\delta_{i,j}$ is the Kronecker symbol, which equals 1 for $i = j$, and is zero otherwise. The modal equations (21) are uncoupled. Applying the central difference scheme for the time integration yields the following discrete representation [9, p. 400].

$$\frac{1}{\Delta t^2}\tilde{x}^{n+1} + \left(\mu_j + 2\frac{\xi_j\sqrt{\mu_j}}{\Delta t} - 2\frac{1}{\Delta t^2}\right)\tilde{x}^n + \left(\frac{1}{\Delta t^2} - 2\frac{\xi_j\sqrt{\mu_j}}{\Delta t}\right)\tilde{x}^{n-1} = \mathbf{v}_j^T\mathbf{f}_{\text{ext}}^n \tag{24}$$

The above can be put into the same form as (17), yielding the amplification matrix $\mathbf{A}_j$. Defining $\tilde{\mathbf{y}}^n = \begin{bmatrix}\tilde{x}^n \\ \tilde{x}^{n-1}\end{bmatrix}$,

$$\tilde{\mathbf{y}}^{n+1} = \begin{bmatrix}\tilde{x}^{n+1} \\ \tilde{x}^n\end{bmatrix} = \underbrace{\begin{bmatrix}2 - 2\xi_j\sqrt{\mu_j}\Delta t - \Delta t^2\mu_j & 2\xi_j\sqrt{\mu_j}\Delta t - 1 \\ 1 & 0\end{bmatrix}}_{\mathbf{A}_j}\begin{bmatrix}\tilde{x}^n \\ \tilde{x}^{n-1}\end{bmatrix} + \begin{bmatrix}\mathbf{v}_j^{n+1^T}\mathbf{f}_{\text{ext}}^n \\ 0\end{bmatrix}. \tag{25}$$

As only solutions with $\Delta t > 0$ are relevant, the critical value $\|\mathbf{A}_j\|_2 = \max\left(\left|\lambda_1^{A_j}\right|, \left|\lambda_2^{A_j}\right|\right) = 1$ is obtained for [9, p. 401]

$$\Delta t_{\text{crit},j} = \frac{2}{\sqrt{\mu_j}}\left(\sqrt{\xi_j^2 + 1} - \xi_j\right) \tag{26}$$

and any values $0 < \Delta t < \min_j \Delta t_{\text{crit},j}$ will fulfil the stability conditions for *all* modal equations. A more in-depth step-by-step derivation of the above condition can be obtained by applying the $z$-transform and the Routh-Hurwitz criterion [9, 15, 16].

Notice that the computation of $\Delta t_{\text{crit}}$ requires the knowledge of the largest eigenvalue $\mu_m$ of the matrix $\mathbf{M}^{-1}\mathbf{K}$. It is therefore not necessary to compute all eigenvalues, and highly efficient methods can be used instead of a full eigendecomposition. These include Lanczos or Arnoldi iterations, or randomized methods [17]. However, the direct computation of $\mu_m$ would require the assembly of the full stiffness matrix $\mathbf{K}$, which is a costly operation that is usually not performed in most explicit FEM solvers. Instead, it is possible to derive computationally cheap – albeit somewhat conservative – bounds on $\mu_m$ and $\Delta t_{\text{crit}}$ by computing eigenvalues at the element level. The underlying theory is detailed in the next section.

### 3.4. Rayleigh-Ritz quotient and eigenvalue separation property

The Rayleigh quotient of a symmetric matrix $\mathbf{A} \in \mathbb{R}^{m\times m}$ and a nonzero vector $\mathbf{x} \in \mathbb{R}^m$ is defined by [18, pp. 230ff.]

$$R_A(\mathbf{A}, \mathbf{x}) = \frac{\mathbf{x}^T\mathbf{A}\mathbf{x}}{\mathbf{x}^T\mathbf{x}}. \tag{27}$$

If $\mathbf{x}$ is an eigenvector of $\mathbf{A}$, the Rayleigh quotient simplifies to its corresponding eigenvalue. The Rayleigh quotient is bounded as follows (notice that $\mathbf{A}$ is real and symmetric and therefore automatically Hermitian, with only real eigenvalues).

**Theorem.** *Let* $\mathbf{A} \in \mathbb{R}^{m\times m}$ *be a symmetric matrix, and let* $\lambda_i$ *denote the i-th eigenvalue of* $\mathbf{A}$ *in increasing order. Then*

$$\lambda_1 \leq \frac{\mathbf{x}^T\mathbf{A}\mathbf{x}}{\mathbf{x}^T\mathbf{x}} \leq \lambda_m. \tag{28}$$

A proof of this theorem is given in [18, pp. 230ff.]. An important consequence for finite element methods is the Rayleigh nesting theorem. The following version is due to [9]. A proof is given in [19].

**Theorem.** *(Rayleigh-Fried-Belytschko)*
*Let $\lambda_i$ denote the eigenvalues of the generalized eigenvalue problem $\mathbf{Ax} = \lambda\mathbf{Bx}$ in increasing order, where $\mathbf{A} \in \mathbb{R}^{m\times m}$ and $\mathbf{B} \in \mathbb{R}^{m\times m}$ are both symmetric, and constructed from real symmetric matrices $\mathbf{A}_e$ and $\mathbf{B}_e$ by*

$$\mathbf{A} = \sum_e \mathbf{L}_e^T \mathbf{A}_e \mathbf{L}_e \qquad \mathbf{B} = \sum_e \mathbf{L}_e^T \mathbf{B}_e \mathbf{L}_e, \tag{29}$$

*where $\mathbf{L}_e \in \mathbb{R}^{m_e\times m}$ are the element connectivity matrices (cf. [9]). The eigenvalues $\bar{\lambda}_i$ of the system constrained by $\mathbf{g}^T\mathbf{x} = a$, $\mathbf{g} \in \mathbb{R}^m$ are then nested by the eigenvalues of the unconstrained system, i.e.*

$$\lambda_1 \le \bar{\lambda}_1 \le \lambda_2 \le \bar{\lambda}_2 \le \cdots \le \lambda_n. \tag{30}$$

*A consequence of the nesting theorem is the element eigenvalue inequality, which bounds the maximum system eigenvalue by the maximum element eigenvalue. The element eigenvalue inequality can be written as [9, p.403]*

$$|\lambda_{max}| \le |\lambda_{max}^E|, \tag{31}$$

*where $\lambda_{max}$ is the largest eigenvalue of the symmetric generalized eigenvalue problem $\mathbf{Ax} = \lambda\mathbf{Bx}$, and $\lambda_{max}^e = \max_{i,e}\lambda_i^e$ is the largest eigenvalue of any element in the assembly. The element with the largest eigenvalue has the index $e = E$. In this context, we also refer to the work of Fried [20].*

Theorem 3.4 can be applied to analyze the stability of integration methods for the equations of motion (1). The matrices $\mathbf{A}$ and $\mathbf{B}$ in theorem 3.4 then correspond to the system mass and stiffness matrices $\mathbf{M}$ and $\mathbf{K}$. Conversely, $\mathbf{A}_e$ and $\mathbf{B}_e$ correspond to the element mass and stiffness matrices $\mathbf{M}_e$ and $\mathbf{K}_e$. This means that one can bound the largest eigenvalue of the generalized eigenvalue problem $\mathbf{Kx} = \lambda\mathbf{Mx}$, i.e. the largest eigenvalue of $\mathbf{M}^{-1}\mathbf{K}$, by the largest eigenvalue of the element eigenvalue problems $\mathbf{K}_e\mathbf{x}^e = \lambda^e\mathbf{M}_e\mathbf{x}^e$. The advantage is that the element eigenvalues can be calculated (or approximated) at very low computational costs. Moreover, the assembly or storage of the stiffness matrix can be omitted. The critical time step is then approximated using the largest eigenvalue of the critical element $e$. Consider a system consisting of simple 2-node rod element in uniaxial strain, a diagonal mass matrix, and no damping. Combining equation (26) and theorem 3.4, the critical time step for the central difference method can be bounded by [9, p. 404]

$$\Delta t \le \frac{2}{\sqrt{\lambda_{\max}}} \le \frac{2}{\max_{e,i}\lambda_i^e} = \min_{e,i}\frac{l_e}{c_e}, \quad \text{where} \quad \lambda_{\max} = \omega_{\max}^2 \quad \text{and} \quad \lambda_{\max}^e = 4\frac{c_e^2}{l_e^2}. \tag{32}$$

$l_e$ is the length of the rod element $e$, and $c_e$ is the corresponding wave speed, which depends on material parameters such as the density or the Young's modulus. Condition (32) is also referred to as the Courant-Friedrichs-Lewy (CFL) condition [21]. Although the bounds obtained by the eigenvalue inequality are conservative, they are often used in explicit FEM solvers as they avoid the assembly of the full stiffness matrix $\mathbf{K}$, and the computation of even just a single eigenvalue of the system eigenvalue problem $\mathbf{Kx} = \lambda\mathbf{Mx}$ can be too computationally intensive when the system is very large [9, p. 403]. It is also possible to use the even simpler 1D approximation using the edge length and the wave speed of the critical element. An additional safety factor is often included to account for any errors introduced by this approximation, as well as the usage of linear stability theory [10].

## 4. FURTHER EIGENVALUE THEOREMS FOR HERMITIAN MATRICES

This section presents a couple of important eigenvalue theorems for real symmetric (Hermitian) matrices which will become important for the stability analysis of the reduced and hyper reduced systems in sections 5 and 6.

### 4.1. Courant-Fischer theorem and Weyl's inequality

The Courant-Fischer theorem, also commonly referred to as *min-max principle* [22, p. 58], is closely related to theorem 3.4. Proofs of both of the following theorems can be found in [18, pp. 230ff.] and [22, pp. 57ff.], or in the original works of Fischer [23] and Courant [24]. We briefly highlight the most important results for this work.

**Theorem.** *(Magnus)*
*Let* $\mathbf{A} \in \mathbb{R}^{m\times m}$ *be a symmetric matrix, and let* $\lambda_i$ *denote the i-th eigenvalue of* $\mathbf{A}$ *in increasing order. Let further* $\mathbf{B} \in \mathbb{R}^{m\times(k-1)}$ *and* $\mathbf{C} \in \mathbb{R}^{m\times(m-k)}$ *be arbitrary real matrices. Then*

$$\min_{\mathbf{B}^T\mathbf{x}=0} \frac{\mathbf{x}^T\mathbf{A}\mathbf{x}}{\mathbf{x}^T\mathbf{x}} \leq \lambda_k \qquad \textit{and} \qquad \max_{\mathbf{C}^T\mathbf{x}=0} \frac{\mathbf{x}^T\mathbf{A}\mathbf{x}}{\mathbf{x}^T\mathbf{x}} \geq \lambda_k \tag{33}$$

**Theorem.** *(Courant-Fischer)*
*Let* $\mathbf{A} \in \mathbb{R}^{m\times m}$ *be a symmetric matrix, and let* $\lambda_i$ *denote the i-th eigenvalue of* $\mathbf{A}$ *in increasing order. Then, the eigenvalues may also be defined as*

$$\lambda_k = \max_{\mathbf{B}^T\mathbf{B}=\mathbf{I}_{k-1}} \min_{\mathbf{B}^T\mathbf{x}=0} \frac{\mathbf{x}^T\mathbf{A}\mathbf{x}}{\mathbf{x}^T\mathbf{x}} = \min_{\mathbf{C}^T\mathbf{C}=\mathbf{I}_{n-k}} \max_{\mathbf{C}^T\mathbf{x}=0} \frac{\mathbf{x}^T\mathbf{A}\mathbf{x}}{\mathbf{x}^T\mathbf{x}} \tag{34}$$

*where* $\mathbf{B} \in \mathbb{R}^{m\times(k-1)}$ *and* $\mathbf{C} \in \mathbb{R}^{m\times(m-k)}$ *are real matrices.*

An interesting consequence of the Courant-Fischer theorem is the existence of lower and upper bounds on the eigenvalues of the sum of two symmetric real (or Hermitian) matrices. This is summarized in Weyl's inequality theorem [22, pp. 62ff.], initially described in [25].

**Theorem.** *(Weyl)*
*Let* $\mathbf{A} \in \mathbb{R}^{m\times m}$ *and* $\mathbf{B} \in \mathbb{R}^{m\times m}$ *both be symmetric and real-valued (or Hermitian), and let* $\lambda_i^A$ *denote the i-th eigenvalue of* $\mathbf{A}$*, and* $\lambda_i^B$ *the i-th eigenvalue of* $\mathbf{B}$*, in increasing order. Let further* $\lambda_i^{A+B}$ *be the i-th eigenvalue of* $\mathbf{A}+\mathbf{B}$*. Then, for* $1 \leq j \leq m-i+1 \leq m$ *and* $1 \leq l \leq i$*, the following inequalities hold [26]*

$$\lambda_{i-l+1}^A + \lambda_l^B \leq \lambda_i^{A+B} \leq \lambda_{i+j-1}^A + \lambda_{m-j+1}^B, \tag{35}$$

*and thus (for* $1 \leq i \leq m$*)*

$$\lambda_i^A + \lambda_1^B \leq \lambda_i^{A+B} \leq \lambda_i^A + \lambda_m^B \qquad \textit{and} \qquad \lambda_1^A + \lambda_i^B \leq \lambda_i^{A+B} \leq \lambda_m^A + \lambda_i^B. \tag{36}$$

For further bounds of the eigenvalues and singular values of the sum or the product of Hermitian matrices, we refer to [22, chapter III], and [26–29].

### 4.2. Cauchy interlacing theorem and Poincaré separation theorem

**Theorem.** *(Cauchy)*
*Let* $\mathbf{V} \in \mathbb{R}^{m\times(m-1)}$ *be a semi-orthogonal matrix with orthonormal columns. Let further* $\mathbf{A} \in \mathbb{R}^{m\times m}$ *be a symmetric matrix, and let* $\lambda_i$ *denote the i-th eigenvalue of* $\mathbf{A}$ *in increasing order. Denoting the i-th eigenvalue of* $\mathbf{V}^T\mathbf{A}\mathbf{V}$ *by* $\tilde{\lambda}_i$*, again in increasing order, the following chain of inequalities holds.*

$$\lambda_1 \leq \tilde{\lambda}_1 \leq \lambda_2 \leq \tilde{\lambda}_2 \leq \cdots \leq \lambda_{m-1} \leq \tilde{\lambda}_{m-1} \leq \lambda_m \tag{37}$$

This version of the Cauchy interlacing theorem is given in [30, p.72], along with a brief proof. Notice the resemblance to the Rayleigh nesting theorem (theorem 3.4).

**Theorem.** *(Poincaré)*
*Let* $\mathbf{V} \in \mathbb{R}^{m\times k}$*,* $k \leq m$*, be a real semi-orthogonal matrix with orthonormal column vectors, and let* $\mathbf{A} \in \mathbb{R}^{m\times m}$ *be a real symmetric matrix. Then the eigenvalues of the transformed matrix* $\mathbf{V}^T\mathbf{A}\mathbf{V}$ *separate the eigenvalues of* $\mathbf{A}$ *such that*

$$\lambda_i \leq \tilde{\lambda}_i \leq \lambda_{m-k+i}; \quad i \leq k, \tag{38}$$

*where* $\lambda_i$ *denotes the i-th eigenvalue of* $\mathbf{A}$*, and* $\tilde{\lambda}_i$ *is the i-th eigenvalue of* $\mathbf{V}^T\mathbf{A}\mathbf{V}$*.*

A proof of this theorem is given in [18, p. 236]. It is closely related to the Cauchy interlacing theorem. An important corollary for principal submatrices of $\mathbf{A}$ is also given by [18].

**Corollary.** *The Poincaré theorem extends to the case where* $\mathbf{V} = \mathbf{W} \in \mathbb{R}^{m\times m}$ *is an idempotent symmetric* $m \times m$ *matrix of rank* $k \leq m$*. The idempotent property is defined as* $\mathbf{W}\mathbf{W} = \mathbf{W}$*. Let* $\mathbf{A} \in \mathbb{R}^{m\times m}$ *again be a symmetric matrix, with* $\lambda_i$ *being its i-th eigenvalue, in increasing order. Then, the k nonzero eigenvalues* $\tilde{\lambda}_i$ *of the matrix* $\mathbf{W}\mathbf{A}\mathbf{W}$ *satisfy the separation property*

$$\lambda_i \leq \tilde{\lambda}_i \leq \lambda_{m-k+i}. \tag{39}$$

*The remaining* $m-k$ *eigenvalues of* $\mathbf{W}\mathbf{A}\mathbf{W}$ *are zero. A proof of this corollary is given in [18, p. 237].*

**Corollary.** *Let $\mathbf{A} \in \mathrm{I\!R}^{m\times m}$ be a real symmetric matrix, and let $\mathbf{A}_k \in \mathrm{I\!R}^{k\times k}, k \le m$ be a principal submatrix of $\mathbf{A}$. This means that $\mathbf{A}_k$ is obtained by removing any set of $(m-k)$ rows of $\mathbf{A}$, as well as the same $m-k$ columns. Then, denoting the i-th eigenvalue of $\mathbf{A}$ by $\lambda_i$, and the i-th eigenvalue of $\mathbf{A}_k$ by $\tilde{\lambda}_i$,*

$$\lambda_i \le \tilde{\lambda}_i \le \lambda_{m-k+i}; \quad i \le k. \tag{40}$$

*Proof.* In line with [18], we can let $\mathbf{V} \in \mathrm{I\!R}^{m\times k}$ be a truncated identity matrix $\mathbf{V} = \begin{pmatrix} \mathbf{I}_k \\ \mathbf{0} \end{pmatrix}$, or any row permutation thereof. Then $\mathbf{V}^T\mathbf{A}\mathbf{V}$ is a principal submatrix of $\mathbf{A}$, and theorem 4.2 applies. □

The above corollary is a generalized version of the eigenvalue separation property described by Bathe [8, pp. 63ff]. A consequence of this property is that the characteristic polynomials of $\mathbf{A}_k\mathbf{x} = \lambda\mathbf{x}$ represent a Sturm sequence [8, p. 65].

### 4.3. Poincaré separation theorem for singular values

In fact, Rao [27] proved that it is possible to generalize theorem 4.2 for the singular values of a matrix $\mathbf{A}$. The original formulation of the theorem is not restricted to real matrices, and makes use of the conjugate transpose instead of the transpose operator.

**Theorem.** *(Rao)*
*Let $\mathbf{A} \in \mathrm{I\!R}^{m\times n}$. Let further $\mathbf{B} \in \mathrm{I\!R}^{m\times k}$ and $\mathbf{C} \in \mathrm{I\!R}^{n\times r}$ be orthonormal matrices such that $\mathbf{B}^T\mathbf{B} = \mathbf{I}_k$ and $\mathbf{C}^T\mathbf{C} = \mathbf{I}_r$. The i-th singular value of $\mathbf{A}$ is denoted by $\sigma_i$, in increasing order. Conversely, let $\tilde{\sigma}_i$ denote the i-th singular value of $\mathbf{B}^T\mathbf{A}\mathbf{C}$. Then the following separation inequality holds.*

$$\sigma_{m+n-r-k+i} \le \tilde{\sigma}_i \le \sigma_i \qquad \textit{for} \qquad 1 \le i \le \min(r,k). \tag{41}$$

**Corollary.** *Notice that the above theorem also covers the transform of a symmetric real matrix $\mathbf{A} \in \mathrm{I\!R}^{m\times m} \to \mathbf{V}^T\mathbf{A}\mathbf{U}$, where $\mathbf{V} \in \mathrm{I\!R}^{m\times k}$ and $\mathbf{U} \in \mathrm{I\!R}^{m\times k}$ have orthonormal columns such that $\mathbf{V}^T\mathbf{V} = \mathbf{I}_k = \mathbf{U}^T\mathbf{U}$.*

*Proof.* $\mathbf{V}^T\mathbf{A}\mathbf{U}$ is square, but for $\mathbf{U} \ne \mathbf{V}$ it is usually not symmetric and the singular values of $\mathbf{V}^T\mathbf{A}\mathbf{U}$ are in general not identical to its eigenvalues. However, the following are true:

- as $\mathbf{A}$ is symmetric, its singular values $\sigma_i$ are identical to its eigenvalues $\lambda_i$. This is because the singular values of $\mathbf{A}$ are the square roots of the eigenvalues of $\mathbf{A}^T\mathbf{A} = \mathbf{A}\mathbf{A}$. [3]
- The maximum singular value of any arbitrary matrix $\mathbf{B}$ is always larger than or equal to its largest eigenvalue (in absolute value). To show this property, recall that the spectral norm of a matrix is equal to its largest singular value, i.e. $\sigma_m = \sqrt{\|\mathbf{A}^T\mathbf{A}\|_2}$. Then, for any eigenvalue $\lambda_X$ and corresponding eigenvector $\mathbf{v}_X$ of the matrix $\mathbf{X} = \mathbf{V}^T\mathbf{A}\mathbf{U}$,

$$|\lambda_X| = \|\lambda\mathbf{v}_X\|_2 = \|\mathbf{X}\mathbf{v}_X\|_2 = \left(\mathbf{v}_X^T\mathbf{X}^T\mathbf{X}\mathbf{v}_X\right)^{1/2} \le \left\|\mathbf{X}^T\mathbf{X}\right\|_2^{1/2}\left(\mathbf{v}_X^T\mathbf{v}_X\right)^{1/2} = \left\|\mathbf{X}^T\mathbf{X}\right\|_2^{1/2}. \tag{42}$$

A version of this proof was initially given by M. Argerami [31]

Therefore the largest eigenvalue (in absolute magnitude) of $\mathbf{X} = \mathbf{V}^T\mathbf{A}\mathbf{U}$ is always smaller than or equal to its largest singular value. The largest singular value of $\mathbf{V}^T\mathbf{A}\mathbf{U}$ is, in turn, bounded by the largest eigenvalue of $\mathbf{A}$ (theorem 4.3). Hence the largest eigenvalue of $\mathbf{A}$ will always be larger than the largest absolute value of any eigenvalue of $\mathbf{V}^T\mathbf{A}\mathbf{U}$. □

## 5. NUMERICAL STABILITY ANALYSIS OF GALERKIN REDUCED-ORDER MODELS

Having established the relationship of numerical stability to the eigenvalues of the system matrix $\mathbf{M}^{-1}\mathbf{K}$, this section extends the stability analysis to reduced-order models without hyper reduction. It is proved, for the first time that the authors are aware, that a Galerkin projection generally increases the critical time step limit and therefore enables the use of larger time steps within the reduced-order model.

[3] Let $\mathbf{A} = \mathbf{Q}\mathbf{\Lambda}\mathbf{Q}^{-1}$ be the eigendecomposition of $\mathbf{A}$, where the diagonal entries of $\mathbf{\Lambda}$ are the eigenvalues of $\mathbf{A}$. Then, $\mathbf{A}\mathbf{A} = \mathbf{Q}\mathbf{\Lambda}\mathbf{Q}^{-1}\mathbf{Q}\mathbf{\Lambda}\mathbf{Q}^{-1} = \mathbf{Q}\mathbf{\Lambda}^2\mathbf{Q}^{-1}$.

### 5.1. Motivation

As pointed out in the introduction, it has been observed by many authors that reduced-order models using Galerkin projections and proper orthogonal decomposition (POD) tend to have better stability properties for explicit time integration methods than the original full-order model (FOM), enabling the choice of larger time steps to speed up the simulation [1–6]. The use of larger time steps can be associated with significant speed-up factors that would otherwise be impossible without hyper reduction when explicit time integration is used [5]. For the analysis of the dynamic impact of a circular Taylor-type bar, Farhat et al. [5] report a speed-up factor of 118 of the explicitly integrated reduced-order model, compared to the corresponding FOM simulation. This was largely due to much larger critical time step sizes. While this overcomes one of the major drawbacks of explicit integration methods, the reasons for this phenomenon are still largely unknown, and it lacks thorough mathematical analysis. Taylor et al. [4] state that the larger critical time step size is likely due to the fact that the higher frequencies of the accelerations are filtered out by the projection onto the ROB. Farhat et al. [5] and Krysl et al. [2] give explanations with similar reasoning. This section aims to back these explanations with rigorous mathematical analysis, and to provide a better understanding of the numerical stability of projection-based reduced-order methods.

Finally, although not the subject of this paper, we would like to emphasize that similar phenomena can be observed for many problems which are described using more "global" shape functions or modes. One example is Isogeometric Analysis [32], where recent research has found that the interior isogeometric elements of a structure benefit from larger critical time steps due to an increase in the size of the basis function support [33].

### 5.2. Galerkin projections and reduced-order form of the equations of motion

The projection-based ROMs in this work make use of a reduced basis matrix $\mathbf{V} \in \mathbb{R}^{m\times k}$, where $k$ is much smaller than $m$. The orthonormal columns of $\mathbf{V}$ contain "empirical modes" which are often obtained from a truncated singular value decomposition (SVD) of a snapshot matrix [34]. The snapshot compression is equivalently referred to as Proper Orthogonal Decomposition (POD). In structural mechanics, the snapshots correspond to observed deformation states collected from a number of a priori training simulations. The high-dimensional variable $\mathbf{x} \in \mathbb{R}^m$ can then be approximated using the relation

$$\mathbf{x}(t) = \mathbf{x}^0 + \mathbf{V}\tilde{\mathbf{x}}(t), \tag{43}$$

where $\mathbf{x}^0 \in \mathbb{R}^m$ is the vector of initial displacements, and $\tilde{\mathbf{x}} \in \mathbb{R}^k$ is the vector of the reduced DoFs. $\tilde{\mathbf{x}}$ is equivalently referred to as the vector of reduced coordinates, or the projection of $\mathbf{x}$ onto the columns of $\mathbf{V}$. For the sake of simplicity, we restrict the following analysis to global subspace approximations. The extension to local ROBs [35] can be carried out analogously. Without loss of generality, we further assume $\mathbf{x}^0 = \mathbf{0}$, which simplifies (43) to $\mathbf{x}(t) = \mathbf{V}\tilde{\mathbf{x}}(t)$. Inserting (43) into (1) then yields

$$\mathbf{M}(\dot{\mathbf{x}}, \mathbf{x}, t)\ \mathbf{V}\ddot{\tilde{\mathbf{x}}} + \mathbf{C}(\dot{\mathbf{x}}, \mathbf{x}, t)\ \mathbf{V}\dot{\tilde{\mathbf{x}}} + \mathbf{K}(\dot{\mathbf{x}}, \mathbf{x}, t)\ \mathbf{V}\tilde{\mathbf{x}} = \mathbf{f}_{\text{ext}}(\dot{\mathbf{x}}, \mathbf{x}, t). \tag{44}$$

This describes an overdetermined system with $m$ equations and $k \ll m$ DoFs, giving rise to a residual vector $\mathbf{r} \in \mathbb{R}^m$. The subsequent Galerkin projection enforces the residual to be orthogonal to the space spanned by the columns of $\mathbf{V}$ [36] by solving $\mathbf{V}^T\mathbf{r}(\tilde{\mathbf{x}}, \dot{\tilde{\mathbf{x}}}, t) = 0$, i.e.

$$\mathbf{V}^T\mathbf{M}(\dot{\mathbf{x}}, \mathbf{x}, t)\ \mathbf{V}\ddot{\tilde{\mathbf{x}}} + \mathbf{V}^T\mathbf{C}(\dot{\mathbf{x}}, \mathbf{x}, t)\ \mathbf{V}\dot{\tilde{\mathbf{x}}} + \mathbf{V}^T\mathbf{K}(\dot{\mathbf{x}}, \mathbf{x}, t)\ \mathbf{V}\tilde{\mathbf{x}} = \mathbf{V}^T\mathbf{f}_{\text{ext}}(\dot{\mathbf{x}}, \mathbf{x}, t) \tag{45}$$

In the following, we again focus on the linearized form of the equations (section 3.2). The linearized reduced-order equations then have a form which is strikingly similar to (12).

$$\underbrace{\mathbf{V}^T\mathbf{M}\mathbf{V}}_{\mathbf{M}_r}\ddot{\tilde{\mathbf{x}}} + \underbrace{\mathbf{V}^T\mathbf{C}\mathbf{V}}_{\mathbf{C}_r}\dot{\tilde{\mathbf{x}}} + \underbrace{\mathbf{V}^T\mathbf{K}\mathbf{V}}_{\mathbf{K}_r}\tilde{\mathbf{x}} = \underbrace{\mathbf{V}^T\mathbf{f}_{\text{ext}}}_{\mathbf{f}_{\text{ext},r}} \tag{46}$$

As is typical for explicit FEM, we assume a lumped (diagonal) positive mass matrix, as well as symmetric and positive semi-definite damping and stiffness matrices. The columns of $\mathbf{V}$ form an orthonormal basis. This implies $\mathbf{V}^T\mathbf{V} = \mathbf{I}$. The reduced basis $\mathbf{V}$ is also often normalized to $\mathbf{M}$ in practice. In this case the relationship $\mathbf{V}^T\mathbf{M}\mathbf{V} = \mathbf{I}$ holds. We adopt this convention throughout the remainder of this paper, and will furthermore assume a constant mass matrix $\mathbf{M}$ hereafter. $\mathbf{V}_M = \mathbf{M}^{-1/2}\mathbf{V}$ is used to denote the $\mathbf{M}$-orthonormalized version of $\mathbf{V}$ [4]. The reduced equations can then be written as

$$\underbrace{\mathbf{I}}_{\mathbf{M}_r}\ddot{\tilde{\mathbf{x}}} + \underbrace{\mathbf{V}_M^T\mathbf{C}\mathbf{V}_M}_{\mathbf{C}_r}\dot{\tilde{\mathbf{x}}} + \underbrace{\mathbf{V}_M^T\mathbf{K}\mathbf{V}_M}_{\mathbf{K}_r}\tilde{\mathbf{x}} = \underbrace{\mathbf{V}_M^T\mathbf{f}_{\text{ext}}}_{\mathbf{f}_{\text{ext},r}}. \tag{47}$$

We refer to Farhat et al. [5] for details regarding the treatment of rotational DoFs as opposed to displacement DoFs, which is particularly important for models with configuration dependent mass matrices.

### 5.3. Properties of the reduced system matrices

In fact, one can argue that the structure of the equations (46) and (47) is so similar to (12) that one can simply carry out the stability analysis analogously to the FOM. In order to do this, however, it is important to ensure that crucial properties of the system matrices are preserved in the reduced-order equations. These are above all the symmetry and positive definiteness of $\mathbf{M}_r$, $\mathbf{C}_r$ and $\mathbf{K}_r$. It turns out that the above Galerkin projection does indeed preserve these properties, both if $\mathbf{V}$ with $\mathbf{V}^T\mathbf{V} = \mathbf{I}$, or if $\mathbf{V}_M$ with $\mathbf{V}_M^T\mathbf{M}\mathbf{V}_M = \mathbf{I}$ are used as a reduced basis. The proofs are presented in the following.

**Theorem.** *(Symmetry of $\mathbf{V}^T\mathbf{A}\mathbf{V}$)*
*Let $\mathbf{V} \in \mathbb{R}^{m\times k}, k \leq m$ be a real semi-orthogonal matrix with orthonormal column vectors, and let $\mathbf{A} \in \mathbb{R}^{m\times m}$ be a real symmetric matrix. Then, $\mathbf{V}^T\mathbf{A}\mathbf{V}$ is also a real symmetric matrix, and the transform $\mathbf{A} \to \mathbf{V}^T\mathbf{A}\mathbf{V}$ thus preserves symmetry.*

*Proof.* Let $\mathbf{B} = \mathbf{V}^T\mathbf{A}\mathbf{V}$. Then the entries of $\mathbf{B}$ are given by

$$B_{i,j} = \sum_{p=1}^{m} \left(\mathbf{V}^T\mathbf{A}\right)_{i,p} V_{p,j} = \sum_{p=1}^{m} \left(\sum_{q=1}^{m} V_{q,i}A_{q,p}\right) V_{p,j} = \sum_{p=1}^{m}\sum_{q=1}^{m} V_{p,j}V_{q,i}A_{q,p}. \tag{48}$$

Conversely, the entries $B_{j,i}$ are given by (using the associative property of matrix multiplication)

$$B_{j,i} = \sum_{p=1}^{m} V_{p,j}\left(\mathbf{A}\mathbf{V}\right)_{p,i} = \sum_{p=1}^{m} V_{p,j}\left(\sum_{q=1}^{m} A_{p,q}V_{q,i}\right) = \sum_{p=1}^{m}\sum_{q=1}^{m} V_{p,j}V_{q,i}A_{p,q}. \tag{49}$$

Because of the symmetry of $\mathbf{A}$, we have $A_{p,q} = A_{q,p}$, and therefore $B_{i,j} = B_{j,i}$. $\mathbf{B} = \mathbf{V}^T\mathbf{A}\mathbf{V}$ is hence symmetric. □

**Corollary.** *Theorem 5.3 also holds for any other real (or complex) $m \times k$ matrix $\mathbf{V}$, under otherwise identical assumptions. This is true because no particular property of $\mathbf{V}$ is exploited in the proof of theorem 5.3. The associative property of matrix multiplication and the symmetry of $\mathbf{A}$ are sufficient for the proof. In particular, it follows that theorem 5.3 is also true if $\mathbf{V}_M \in \mathbb{R}^{m\times k}$ is used instead of $\mathbf{V}$, where $\mathbf{V}_M$ is an $\mathbf{M}$-orthonormal matrix (with $\mathbf{M}$ diagonal and positive definite), i.e. $\mathbf{V}_M^T\mathbf{M}\mathbf{V}_M = \mathbf{I}_{k\times k}$.*

**Theorem.** *(Real eigenvalues)*
*Let $\mathbf{V} \in \mathbb{R}^{m\times k}, k \leq m$ be a real semi-orthogonal matrix with orthonormal column vectors, and let $\mathbf{A} \in \mathbb{R}^{m\times m}$ be a real symmetric matrix. Then, the eigenvalues of $\mathbf{V}^T\mathbf{A}\mathbf{V}$ are all real.*

*Proof.* Since $\mathbf{A}$ is symmetric, theorem 5.3 ensures that $\mathbf{V}^T\mathbf{A}\mathbf{V}$ is also symmetric. As every symmetric matrix is also a Hermitian matrix, and Hermitian matrices only have real eigenvalues, the eigenvalues of $\mathbf{A}$ and $\mathbf{V}^T\mathbf{A}\mathbf{V}$ must also be real. □

**Theorem.** *(Extension of the Poincaré separation theorem to $\mathbf{M}$-orthonormal matrices)*
*The Poincaré separation theorem (theorem 4.2) can be extended for the case where an $\mathbf{M}$-orthonormal matrix $\mathbf{V}_M \in \mathbb{R}^{m\times k}$ is used instead of $\mathbf{V}$, with $\mathbf{M}$ positive diagonal and $\mathbf{V}_M^T\mathbf{M}\mathbf{V}_M = \mathbf{I}_{k\times k}$, under otherwise identical assumptions. The eigenvalues of $\mathbf{V}_M^T\mathbf{A}\mathbf{V}_M$ then separate the eigenvalues of the matrix $\mathbf{M}^{-1}\mathbf{A}$, such that*

$$\lambda_i \leq \tilde{\lambda}_i \leq \lambda_{m-k+i}, \tag{50}$$

*where $\mathbf{A} \in \mathbb{R}^{m\times m}$ is symmetric, $\lambda_i$ is the i-th eigenvalue of $\mathbf{M}^{-1}\mathbf{A}$, and $\tilde{\lambda}_i$ is the i-th eigenvalue of $\mathbf{V}_M^T\mathbf{A}\mathbf{V}_M$, in increasing order.*

---

[4] For every $\mathbf{M}$-orthonormal basis $\mathbf{V}_M$, there is an orthonormal matrix $\mathbf{V} = \mathbf{M}^{-1/2}\mathbf{V}_M$, as $\mathbf{I} = \mathbf{V}_M^T\mathbf{M}\mathbf{V}_M = \mathbf{V}^T\mathbf{M}^{-1/2}\mathbf{M}\mathbf{M}^{-1/2}\mathbf{V} = \mathbf{V}^T\mathbf{V}$.

*Proof.* To prove this point, let us denote $\mathbf{V}$ as the "standard form" of a semi-orthogonal matrix with $\mathbf{V}^T\mathbf{V} = \mathbf{I}_{k\times k}$. Conversely, $\mathbf{V}_M$ is used to denote its corresponding $\mathbf{M}$-orthonormalized form with $\mathbf{V}_M^T\mathbf{M}\mathbf{V}_M = \mathbf{I}_{k\times k}$. Here, $\mathbf{M}$ is a diagonal and positive definite mass matrix. Therefore, the entries of $\mathbf{V}_M$ must satisfy

$$\mathbf{V}_M^T\mathbf{M}\mathbf{V}_M = \mathbf{V}_M^T\mathbf{M}^{1/2}\mathbf{M}^{1/2}\mathbf{V}_M = (\mathbf{M}^{1/2}\mathbf{V}_M)^T\mathbf{M}^{1/2}\mathbf{V}_M = \mathbf{I}, \tag{51}$$

which is fulfilled for $\mathbf{V}_M = \mathbf{M}^{-1/2}\mathbf{V}$. The matrix $\mathbf{V}_M^T\mathbf{A}\mathbf{V}_M$, with $\mathbf{A}$ real and symmetric, can therefore be written as

$$\mathbf{V}_M^T\mathbf{A}\mathbf{V}_M = (\mathbf{M}^{-1/2}\mathbf{V})^T\mathbf{A}\mathbf{M}^{-1/2}\mathbf{V} = \mathbf{V}^T\mathbf{M}^{-1/2}\mathbf{A}\mathbf{M}^{-1/2}\mathbf{V} \tag{52}$$

Now $\mathbf{M}^{-1/2}\mathbf{A}\mathbf{M}^{-1/2}$ is symmetric, and its eigenvalues are identical to the ones of $\mathbf{M}^{-1}\mathbf{A}$. This can be shown by as follows: First of all, consider two arbitrary square matrices $\mathbf{B} \in \mathbb{R}^{m\times m}$ and $\mathbf{C} \in \mathbb{R}^{m\times m}$. Then, the matrix $\mathbf{BC}$ has the same eigenvalues as the matrix $\mathbf{CB}$, although in general $\mathbf{BC} \neq \mathbf{CB}$. Denoting the $i$-th eigenvalue of $\mathbf{BC}$ with $\theta_i$, and the corresponding $i$-th eigenvector of $\mathbf{BC}$ with $\mathbf{v}_i$, this follows readily from the eigenproblem

$$\theta_i\mathbf{v}_i = \mathbf{BC}\mathbf{v}_i, \tag{53}$$

where we at first consider the case $\theta_i \neq 0$. Pre-multiplying this equation with $\mathbf{C}$ yields

$$\theta_i\mathbf{C}\mathbf{v}_i = \mathbf{CBC}\mathbf{v}_i, \tag{54}$$

i.e. the vector $\mathbf{C}\mathbf{v}_i$ is an eigenvector for the *same* eigenvalue $\theta_i$. This means that $\theta_i$ is an eigenvalue of both $\mathbf{BC}$ and $\mathbf{CB}$. In the case of $\theta_i = 0$, we have

$$\det(\mathbf{BC}) = 0 = \det(\mathbf{B})\det(\mathbf{C}) = \det(\mathbf{CB}). \tag{55}$$

Hence 0 must also be an eigenvalue of $\mathbf{CB}$. Since the above is true for all eigenvalues (all $i$), the eigenvalues of $\mathbf{BC}$ and $\mathbf{CB}$ must be identical. This important result has been described in [37, pp. 41-43].
Now let $\mu_i$ be the $i$-th eigenvalue of $\mathbf{M}^{-1}\mathbf{A}$, and let $\tilde{\mu}_i$ be the $i$-th eigenvalue of $\mathbf{M}^{-1/2}\mathbf{A}\mathbf{M}^{-1/2}$, in increasing order. The matrix $\mathbf{M}^{-1/2}\mathbf{A}\mathbf{M}^{-1/2}$ can be equivalently be written in the form $(\mathbf{M}^{-1/2}\mathbf{A})\mathbf{M}^{-1/2}$. Applying the above result then means that the eigenvalues of $(\mathbf{M}^{-1/2}\mathbf{A})\mathbf{M}^{-1/2}$ are identical to the eigenvalues of $\mathbf{M}^{-1/2}(\mathbf{M}^{-1/2}\mathbf{A}) = \mathbf{M}^{-1}\mathbf{A}$. This automatically implies that the eigenvalues of $\mathbf{M}^{-1/2}\mathbf{A}\mathbf{M}^{-1/2}$ and $\mathbf{M}^{-1}\mathbf{A}$ are the same, and thus $\mu_i = \tilde{\mu}_i$. It can be shown that the eigenvalues are the same as those of $\mathbf{A}\mathbf{M}^{-1}$ by the same reasoning.
Let $\mathbf{G} = \mathbf{M}^{-1/2}\mathbf{A}\mathbf{M}^{-1/2}$. Equation (52) then reads

$$\mathbf{V}_M^T\mathbf{A}\mathbf{V}^M = (\mathbf{M}^{-1/2}\mathbf{V})^T\mathbf{A}\mathbf{M}^{-1/2}\mathbf{V} = \mathbf{V}^T\mathbf{M}^{-1/2}\mathbf{A}\mathbf{M}^{-1/2}\mathbf{V} = \mathbf{V}^T\mathbf{G}\mathbf{V}, \tag{56}$$

where $\mathbf{G}$ is a symmetric matrix (corollary to theorem 5.3) with the same eigenvalue spectrum as $\mathbf{M}^{-1}\mathbf{A}$. As $\mathbf{V}$ is semi-orthogonal, the Poincaré separation theorem (theorem 4.2) now ensures that the eigenvalues of $\mathbf{V}^T\mathbf{G}\mathbf{V}$ separate the eigenvalues of $\mathbf{G}$, and therefore also the eigenvalues of $\mathbf{M}^{-1}\mathbf{A}$:

$$\mu_i \leq \tilde{\mu}_i^* \leq \mu_{m-k+i}, \tag{57}$$

where $\mu_i$ is the $i$-th eigenvalue of $\mathbf{M}^{-1}\mathbf{A}$, and $\tilde{\mu}_i^*$ is the $i$-th eigenvalue of $\mathbf{V}_M^T\mathbf{A}\mathbf{V}_M$, in increasing order. □

**Theorem.** *(Positive-semidefiniteness)*
*Let* $\mathbf{V} \in \mathbb{R}^{m\times k}, k \leq m$, *be a real semi-orthogonal matrix with orthonormal column vectors, and let* $\mathbf{A} \in \mathbb{R}^{m\times m}$ *be a symmetric positive semi-definite matrix. Then,* $\mathbf{V}^T\mathbf{A}\mathbf{V}$ *is also a symmetric positive semi-definite matrix. In addition to the symmetry-preserving property proved in theorem 5.3, the transform* $\mathbf{A} \to \mathbf{V}^T\mathbf{A}\mathbf{V}$ *thus also preserves positive semi-definiteness.*

*Proof.* By theorem 5.3, all eigenvalues of $\mathbf{A}$ and $\mathbf{A} \to \mathbf{V}^T\mathbf{A}\mathbf{V}$ are real-valued. The eigenvalues $\lambda_i$ of $\mathbf{A}$ further satisfy $\lambda_i \geq 0$ due to the positive semi-definite property of $\mathbf{A}$. Let us now denote $\lambda_i$ in increasing order, with $\lambda_m$ as the largest and $\lambda_1$ as the smallest eigenvalue of $\mathbf{A}$. The $k$ eigenvalues of $\mathbf{A} \to \mathbf{V}^T\mathbf{A}\mathbf{V}$ are analogously denoted in increasing order as $\tilde{\lambda}_i$.
The Poincaré separation theorem (theorem 4.2) then yields $\lambda_m \geq \tilde{\lambda}_k$, and in particular $\tilde{\lambda}_1 \geq \lambda_1$. This means that

$$\min_i \tilde{\lambda}_i \geq \min_i \lambda_i \geq 0, \tag{58}$$

and therefore all eigenvalues of $\mathbf{V}^T\mathbf{A}\mathbf{V}$ are also greater than or equal to zero. This immediately implies positive-semidefiniteness of the matrix $\mathbf{V}^T\mathbf{A}\mathbf{V}$. □

**Corollary.** *Analogous to theorem 5.3, the transform* $\mathbf{A} \to \mathbf{V}^T\mathbf{A}\mathbf{V}$ *also preserves positive definiteness if* $\mathbf{A}$ *is symmetric and positive definite, under otherwise identical assumptions.*

**Corollary.** *If* $\mathbf{V}_M \in \mathbb{R}^{m\times k}$ *is an* $\mathbf{M}$*-orthonormal matrix with* $\mathbf{M}$ *positive and diagonal,* $\mathbf{V}_M^T\mathbf{M}\mathbf{V}_M = \mathbf{I}_{k\times k}$*, and otherwise the same assumptions as above, then the matrix* $\mathbf{V}_M^T\mathbf{A}\mathbf{V}_M$ *is positive (semi-)definite if* $\mathbf{M}^{-1}\mathbf{A}$ *is positive (semi-) definite.*

*Proof.* In line with the proof of theorem 5.3, let us again denote $\mathbf{V}$ as the "standard form" of a semi-orthogonal matrix with $\mathbf{V}^T\mathbf{V} = \mathbf{I}_{k\times k}$. $\mathbf{V}_M = \mathbf{M}^{-1/2}$ is used to denote its corresponding $\mathbf{M}$-orthonormalized form with $\mathbf{V}_M^T\mathbf{M}\mathbf{V}_M = \mathbf{I}_{k\times k}$. Then, $\mathbf{V}_M^T\mathbf{A}\mathbf{V}_M$ can be written as (cf. equation (52))

$$\mathbf{V}_M^T\mathbf{A}\mathbf{V}_M = \mathbf{V}^T\mathbf{M}^{-1/2}\mathbf{A}\mathbf{M}^{-1/2}\mathbf{V}.$$

As shown within the proof of theorem 5.3, the eigenvalues of $\mathbf{M}^{-1/2}\mathbf{A}\mathbf{M}^{-1/2}$ are identical to those of $\mathbf{M}^{-1}\mathbf{A}$. If $\mathbf{M}^{-1}\mathbf{A}$ is positive (semi-) definite, theorem 5.3 and the subsequent corollary then ensure positive (semi-) definiteness of $\mathbf{V}_M^T\mathbf{A}\mathbf{V}_M$. □

Theorem 5.3 and the above corollaries also follow from Sylvester's Law of Inertia [38], as in [39, p. 448].

**Corollary.** *Notice that* $\mathbf{M}^{-1}\mathbf{A}$ *is also positive (semi-) definite, as it is the product of two Hermitian matrices, both of which are positive (semi-) definite and thus have eigenvalues* $> 0$ $(\geq 0)$*. Then* $\mathbf{M}^{-1}\mathbf{A}$ *has real eigenvalues and the number of positive (non-negative) eigenvalues is the same as the number of positive (non-negative) eigenvalues in* $\mathbf{M}$ *or* $\mathbf{K}$*. This observation is formulated and proven in [40, p. 109].*

The above theorems ensure that the eigenvalues of any symmetric matrix $\mathbf{A}$ do not change to the extremes if the matrix is transformed to its reduced counterpart $\mathbf{V}^T\mathbf{A}\mathbf{V}$ or $\mathbf{V}_M^T\mathbf{A}\mathbf{V}_M^T$. The transform preserves symmetry and positive (semi-)definiteness, and important eigenvalue separation properties exist. The theorems were generalized to the case where an $\mathbf{M}$-orthonormal matrix $\mathbf{V}_M$ is used instead of $\mathbf{V}$, where $\mathbf{M}$ is a positive diagonal matrix such as the mass matrix. Returning to equations (12) and (46), these findings have the following important implications:

- symmetry and positive definiteness of $\mathbf{V}^T\mathbf{M}\mathbf{V}$ as well as $\mathbf{V}_M^T\mathbf{M}\mathbf{V}_M$, as $\mathbf{M}$ is symmetric and positive definite.
- symmetry and positive semi-definiteness of $\mathbf{V}^T\mathbf{K}\mathbf{V}$ as well as $\mathbf{V}_M^T\mathbf{K}\mathbf{V}_M$, as $\mathbf{K}$ is symmetric and positive semi-definite. $\mathbf{x}^T\mathbf{K}\mathbf{x} \geq 0$ can be interpreted as the energy required to deform a structure by some vector $\mathbf{x}$ against a force $\mathbf{K}\mathbf{x}$. Negative energies would not be physically sensible. Notice that the positive definiteness of the transformed mass matrix can be similarly interpreted in terms of the energy required to accelerate (deform) the structure against its inertia.
- symmetry, positive semi-definiteness and eigenvalue separation properties for $\mathbf{M}^{-1}\mathbf{K}$ and $\mathbf{M}_r^{-1}\mathbf{K}_r = \mathbf{V}_M^T\mathbf{K}\mathbf{V}_M$, when $\mathbf{V}_M$ is used as a reduced basis. This is due to the fact that the eigenvalues of $\mathbf{M}^{-1}\mathbf{K}$ are the same as those of $\mathbf{M}^{-1/2}\mathbf{K}\mathbf{M}^{-1/2}$. More specifically, the eigenvalue separation property ensures that the largest eigenvalue of $\mathbf{V}^T\mathbf{M}^{-1/2}\mathbf{K}\mathbf{M}^{-1/2}\mathbf{V}$ is no greater than the largest eigenvalue of $\mathbf{M}^{-1}\mathbf{K}$.

### 5.4. Stability time step of the reduced-order system

Having established that the Galerkin projection preserves the above matrix properties, it is possible to carry out the stability analysis analogously to section 3.3. For the central difference scheme, the discrete version of the reduced equations of motion (47) reads

$$\underbrace{\mathbf{V}_M^T\mathbf{M}\mathbf{V}_M}_{=\mathbf{I}=\mathbf{M}_r}\ddot{\mathbf{x}}^n + \underbrace{\mathbf{V}_M^T\mathbf{C}\mathbf{V}_M}_{\mathbf{C}_r}\dot{\mathbf{x}}^{n-1/2} + \underbrace{\mathbf{V}_M^T\mathbf{K}\mathbf{V}_M}_{\mathbf{K}_r}\mathbf{x}^n = \mathbf{V}_M^T\mathbf{f}_{\text{ext}}^n. \tag{59}$$

Here, we have again assumed that the reduced basis $\mathbf{V}_M = \mathbf{M}^{-1/2}\mathbf{V}$ has $\mathbf{M}$-orthonormal columns such that $\mathbf{V}_M^T\mathbf{V}_M = \mathbf{I}$. Further, $\mathbf{M}$ is a diagonal and positive mass matrix, and $\mathbf{V}$ has orthonormal columns. Again assuming Rayleigh damping, $\mathbf{C}$ can be written as $\mathbf{C} = a_1\mathbf{M}+a_2\mathbf{K}$, with $a_1 \geq 0, a_2 \geq 0$. Then the transformation $\mathbf{C} \to \mathbf{V}_M^T\mathbf{C}\mathbf{V}_M$ preserves the coefficients $a_1, a_2$, as

$$\mathbf{V}_M^T\mathbf{C}\mathbf{V}_M = \mathbf{V}_M^T\left(a_1\mathbf{M} + a_2\mathbf{K}\right)\mathbf{V}_M = a_1\mathbf{V}_M^T\mathbf{M}\mathbf{V}_M + a_2\mathbf{V}_M^T\mathbf{K}\mathbf{V}_M = a_1\mathbf{I} + a_2\mathbf{V}_M^T\mathbf{K}\mathbf{V}_M. \tag{60}$$

However, the Rayleigh damping coefficients $\xi_j$ of the FOM and $\tilde{\xi}_j$ of the ROM will in general be different, since

$$\xi_j = \frac{a_1}{2\sqrt{\mu_j^{FOM}}} + \frac{a_2\sqrt{\mu_j^{FOM}}}{2}; \qquad \tilde{\xi}_j = \frac{a_1}{2\sqrt{\tilde{\mu}_j^{ROM}}} + \frac{a_2\sqrt{\tilde{\mu}_j^{ROM}}}{2}, \tag{61}$$

where $\mu_j^{FOM}$ is the $j$-th eigenvalue of $\mathbf{M}^{-1}\mathbf{K}$, and $\tilde{\mu}_j^{ROM}$ is the $j$-th eigenvalue of $\mathbf{M}_r^{-1}\mathbf{K}_r = \mathbf{V}_M^T\mathbf{K}\mathbf{V}_M$, which is equal to the $j$-th eigenvalue of $\mathbf{V}^T\mathbf{M}^{-1/2}\mathbf{K}\mathbf{M}^{-1/2}\mathbf{V}$. By the same analysis as in section 3.3, the critical time step of the reduced-order system can be obtained as

$$\Delta t_{\text{crit},j} = \frac{2}{\sqrt{\tilde{\mu}_j^{ROM}}}\left(\sqrt{\tilde{\xi}_j^2 + 1} - \tilde{\xi}_j\right). \tag{62}$$

Due to the Poincaré separation theorem and the eigenvalue separation property of $\mathbf{V}^T\mathbf{M}^{-1/2}\mathbf{K}\mathbf{M}^{-1/2}\mathbf{V}$ with respect to $\mathbf{M}^{-1}\mathbf{K}$, each eigenvalue of $\mathbf{V}^T\mathbf{M}^{-1/2}\mathbf{K}\mathbf{M}^{-1/2}\mathbf{V}$ (now in *decreasing* order) is smaller than or equal to the corresponding eigenvalue of $\mathbf{M}^{-1}\mathbf{K}$, i.e. $\mu_j^{FOM} \geq \tilde{\mu}_j^{ROM}, 1 \leq j \leq k$. In the absence of damping ($\tilde{\xi}_j = 0; \quad 1 \leq j \leq k$), it is therefore obvious that the critical time step size of the ROM cannot exceed the critical time step size of the FOM. The critical time step is furthermore associated with $\mu_1^{FOM}$ ($j$ in decreasing order) of $\mathbf{M}^{-1}\mathbf{K}$, or with $\tilde{\mu}_1^{ROM}$ of $\mathbf{V}^T\mathbf{M}^{-1/2}\mathbf{K}\mathbf{M}^{-1/2}\mathbf{V}$, so that computationally efficient methods for computing only the largest eigenvalue can be used.

The analysis is more involved if damping is present. The minimum critical time step $\Delta t_{\text{crit}} = \min_j \frac{2}{\sqrt{\tilde{\mu}_j^{ROM}}}\left(\sqrt{\tilde{\xi}_j^2 + 1} - \tilde{\xi}_j\right)$ can then be written as a function of $\tilde{\mu}_j^{ROM}$, $a_1$ and $a_2$:

$$\Delta t_{\text{crit}} = \min_j \frac{2}{\sqrt{\tilde{\mu}_j^{ROM}}}\left(\sqrt{\left(\frac{a_1}{2\sqrt{\tilde{\mu}_j^{ROM}}} + \frac{a_2\sqrt{\tilde{\mu}_j^{ROM}}}{2}\right)^2 + 1} - \frac{a_1}{2\sqrt{\tilde{\mu}_j^{ROM}}} - \frac{a_2\sqrt{\tilde{\mu}_j^{ROM}}}{2}\right) \tag{63}$$

In practice, it is important to know for which $j$ the value of $\Delta t_{\text{crit},j}$ is minimal. If $j = 1$, then the critical time step is associated with the largest eigenvalue $\tilde{\mu}_1^{ROM}$ of $\mathbf{V}^T\mathbf{M}^{-1/2}\mathbf{K}\mathbf{M}^{-1/2}\mathbf{V}$. Introducing $x = \sqrt{\tilde{\mu}_j^{ROM}}$, we thus seek the minimum of the function

$$g(x) = \frac{2}{x}\left(\sqrt{\left(\frac{a_1}{2x} + \frac{a_2 x}{2}\right)^2 + 1} - \frac{a_1}{2x} - \frac{a_2 x}{2}\right) = \frac{2}{x}\sqrt{\left(\frac{a_1}{2x} + \frac{a_2 x}{2}\right)^2 + 1} - \frac{a_1}{x^2} - a_2. \tag{64}$$

The influence of damping causes $g(x)$ – and thus the critical time step – to drop below the corresponding curve of the undamped problem. Unfortunately $g(x)$ exhibits a large amount of numerical noise and can only be accurately visualized using extended precision arithmetic (figure 1). The derivative $\frac{d}{dx}g(x) = g'(x)$ for $x \geq 0, a_1 \geq 0, a_2 \geq 0$ is given by

$$g'(x) = \frac{2}{x}\left(\frac{\left(\frac{a_2}{2} - \frac{a_1}{2x^2}\right)\left(\frac{a_1}{2x} + \frac{a_2 x}{2}\right)}{\sqrt{\left(\frac{a_1}{2x} + \frac{a_2 x}{2}\right)^2 + 1}} - \frac{1}{x}\sqrt{\left(\frac{a_1}{2x} + \frac{a_2 x}{2}\right)^2 + 1} + \frac{a_1}{x^2}\right) \tag{65}$$

$$\iff g'(x) = \underbrace{\frac{2}{x^3\sqrt{\left(\frac{a_1}{2x} + \frac{a_2 x}{2}\right)^2 + 1}}}_{>0}\left(\underbrace{-\frac{a_1^2}{2x} - \frac{a_1 a_2 x}{2} - x}_{\leq 0} \underbrace{+ a_1\sqrt{\left(\frac{a_1}{2x} + \frac{a_2 x}{2}\right)^2 + 1}}_{\geq 0}\right). \tag{66}$$

We now prove that $g'(x) \leq 0 \quad \forall x$, and thus the function $g(x)$ is monotonically decreasing. In order for $g'(x)$ to be negative, the factor in brackets needs to become negative, i.e.

$$-\frac{a_1^2}{2x} - \frac{a_1 a_2 x}{2} - x + a_1\sqrt{\left(\frac{a_1}{2x} + \frac{a_2 x}{2}\right)^2 + 1} \leq 0. \tag{67}$$

Rearranging yields

$$a_1\sqrt{\left(\frac{a_1}{2x} + \frac{a_2 x}{2}\right)^2 + 1} \leq \frac{a_1^2}{2x} + \frac{a_1 a_2 x}{2} + x. \tag{68}$$

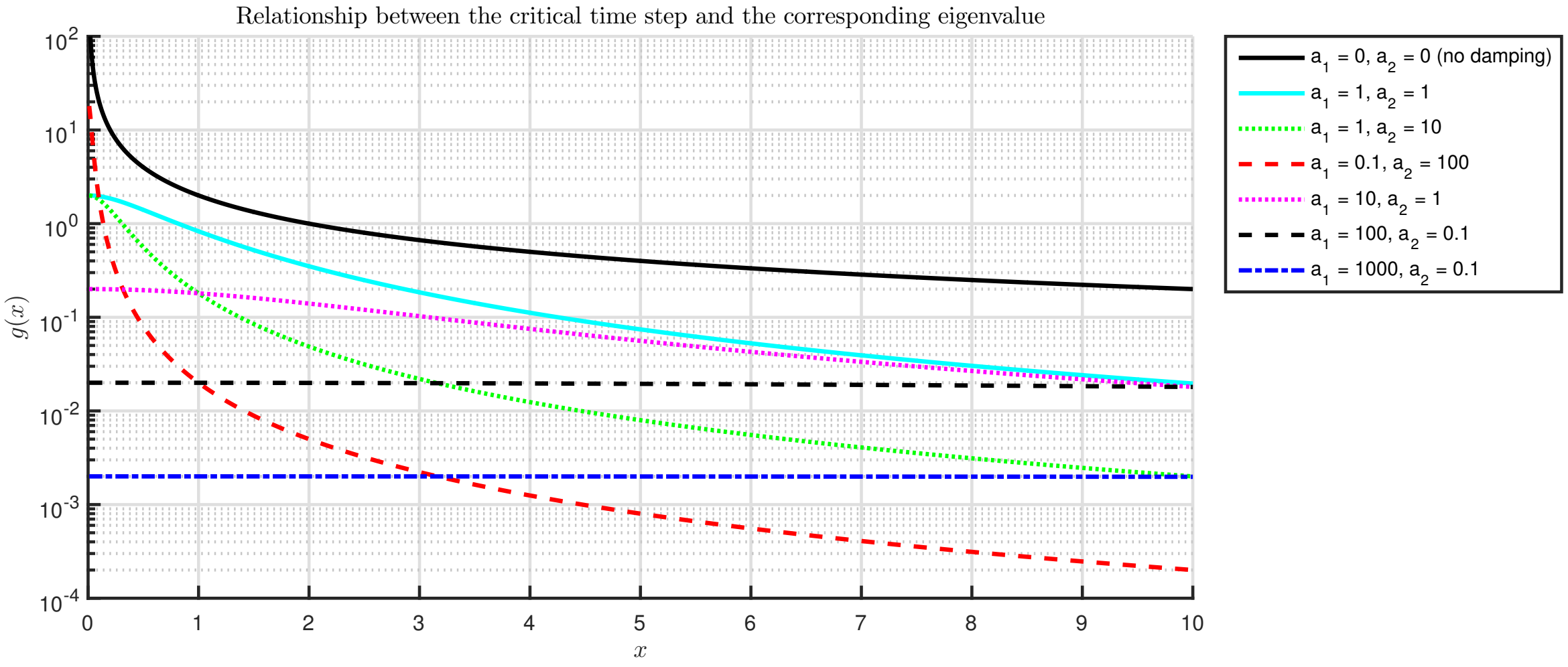

FIG. 1: Critical time step of the ROM for the $j$-th mode, depending on $x = \sqrt{\mu_j^{ROM}}$, plotted using Matlab.

Taking the square $(\cdot)^2$ on both sides, we obtain

$$a_1^2 + \frac{a_1^4}{4x^2} + \frac{a_1^2 a_2^2 x^2}{4} + \frac{a_1^3 a_2}{2} \leq \frac{a_1^4}{4x^2} + \frac{a_1^3 a_2}{4} + \frac{a_1^2}{2} + \frac{a_1^3 a_2}{4} + \frac{a_1^2 a_2^2 x^2}{4} + \frac{a_1 a_2 x^2}{2} + \frac{a_1^2}{2} + \frac{a_1 a_2 x^2}{2} + x^2 \tag{69}$$

$$\iff 0 \leq a_1 a_2 x^2 + x^2, \tag{70}$$

which is always fulfilled since $x \geq 0, a_1 \geq 0, a_2 \geq 0$. Hence $g(x)$ is monotonically decreasing. The limit for $x \to \infty$ is

$$\lim_{x\to\infty} g(x) = \lim_{x\to\infty} \frac{2}{x}\sqrt{\left(\frac{a_1}{2x} + \frac{a_2 x}{2}\right)^2 + 1} - \frac{a_1}{x^2} - a_2 = \lim_{x\to\infty} \sqrt{\frac{a_1^2}{x^4} + 2\frac{a_1 a_2}{x^2} + a_2^2 + \frac{1}{x^2}} - a_2 = 0 \tag{71}$$

The limit for $x \to 0, x \geq 0$ is

$$\lim_{x\to 0} g(x) = \frac{2}{a_1}. \tag{72}$$

As $g(x)$ is monotonically decreasing, it suffices to compute the largest eigenvalue of the system in order to compute the critical time step even if damping is present. $\Delta t_{\text{crit}}$ is given by

$$\Delta t_{\text{crit}}^{ROM} = \min_j \frac{2}{\sqrt{\tilde{\mu}_j^{ROM}}}\left(\sqrt{\tilde{\xi}_j^2 + 1} - \tilde{\xi}_j\right) = \frac{2}{\sqrt{\tilde{\mu}_{\max}^{ROM}}}\left(\sqrt{\left(\frac{a_1}{2\sqrt{\tilde{\mu}_{\max}^{ROM}}} + \frac{a_2\sqrt{\tilde{\mu}_{\max}^{ROM}}}{2}\right)^2 + 1} - \frac{a_1}{2\sqrt{\tilde{\mu}_{\max}^{ROM}}} - \frac{a_2\sqrt{\tilde{\mu}_{\max}^{ROM}}}{2}\right), \tag{73}$$

where $\tilde{\mu}_{\max}^{ROM}$ is the largest eigenvalue of $\mathbf{V}^T\mathbf{M}^{-1/2}\mathbf{K}\mathbf{M}^{-1/2}\mathbf{V}$. Conversely, the critical time step of the FOM can be computed as

$$\Delta t_{\text{crit}}^{FOM} = \min_j \frac{2}{\sqrt{\mu_j^{FOM}}}\left(\sqrt{\xi_j^2 + 1} - \xi_j\right) = \frac{2}{\sqrt{\mu_{\max}^{FOM}}}\left(\sqrt{\left(\frac{a_1}{2\sqrt{\mu_{\max}^{FOM}}} + \frac{a_2\sqrt{\mu_{\max}^{FOM}}}{2}\right)^2 + 1} - \frac{a_1}{2\sqrt{\mu_{\max}^{FOM}}} - \frac{a_2\sqrt{\mu_{\max}^{FOM}}}{2}\right), \tag{74}$$

where $\mu_{\max}^{FOM}$ is the largest eigenvalue of $\mathbf{M}^{-1}\mathbf{K}$. Due to the monotonicity property, larger eigenvalues require smaller time steps. The Poincaré separation theorem ensures, as pointed out in the previous sections, that the largest eigenvalue belonging to the ROM matrix $\mathbf{V}^T\mathbf{M}^{-1/2}\mathbf{K}\mathbf{M}^{-1/2}\mathbf{V}$ cannot be any larger than the largest eigenvalue of the FOM matrix $\mathbf{M}^{-1}\mathbf{K}$. Consequently, the critical stability time step of the ROM can also be no smaller than that of the FOM.

### 5.5. Galerkin reduced-order model stability theorem

The findings of the previous sections can be summarized in a central theorem detailing the computation of the critical stability time step of Galerkin reduced-order models.

**Theorem.** *(Central difference method stability time step for Galerkin ROMs)*
*Let the FOM be a linearized system of the form (46), where $\mathbf{M} \in \mathbb{R}^{m\times m}$ is positive and diagonal, $\mathbf{K} \in \mathbb{R}^{m\times m}$ is positive semi-definite and symmetric, and $\mathbf{C} = a_1\mathbf{M} + a_2\mathbf{K}$ is also symmetric and positive semi-definite. Let further $\mathbf{V}_M \in \mathbb{R}^{m\times k}, k < m$ be a reduced basis with $\mathbf{M}$-orthonormal columns, such that $\mathbf{V}_M^T\mathbf{M}\mathbf{V}_M = \mathbf{I}$, $\mathbf{V}_M = \mathbf{M}^{-1/2}\mathbf{V}$, and $\mathbf{V}^T\mathbf{V} = \mathbf{I}$. Let the ROM be the Galerkin reduced-order model corresponding to the FOM, with reduced basis $\mathbf{V}_M$. The equations of motion of the linearized Galerkin ROM are then given by (47). Then, applying the explicit central difference scheme (section 2.1) to integrate both the FOM and the ROM,*

- *the critical time step of the ROM $\Delta t_{crit}^{ROM}$ cannot be any smaller than the critical stability time step of the FOM $\Delta t_{crit}^{FOM}$*

$$\Delta t_{crit}^{ROM} \geq \Delta t_{crit}^{FOM}. \tag{75}$$

- *the critical time step of the ROM is given by*

$$\Delta t_{crit,k} = \frac{2}{\sqrt{\tilde{\mu}_k^{ROM}}}\left(\sqrt{\tilde{\xi}_k^2 + 1} - \tilde{\xi}_k\right), \qquad \tilde{\xi}_k = \frac{a_1}{2\sqrt{\tilde{\mu}_k^{ROM}}} + \frac{a_2\sqrt{\tilde{\mu}_k^{ROM}}}{2}, \tag{76}$$

*where $\tilde{\mu}_k^{ROM}$ is the largest eigenvalue of the matrix $\mathbf{V}^T\mathbf{M}^{-1/2}\mathbf{K}\mathbf{M}^{-1/2}\mathbf{V}$.*

*The above statements hold regardless of the exact modal content of $\mathbf{V}$, and also regardless of whether $\mathbf{V}$ was obtained from an SVD of displacements, accelerations, or other snapshots.*

Numerical experiments indicate that the above theorem can possibly be extended to damping matrices other than classical Rayleigh damping.

## 6. STABILITY TIME STEP OF HYPER REDUCED MODELS

The previous sections have established that Galerkin reduced-order models with an $\mathbf{M}$-orthonormal reduced basis $\mathbf{V}_M$ do not compromise numerical stability of the equations of motion. Instead, by the eigenvalue separation property, the eigenvalues of the reduced systems usually allow for larger stability time step sizes than in the full model. However, this property alone is often not satisfactory enough to produce large speed-ups. An additional hyper reduction procedure is therefore typically employed to reduce the number of nonlinear function evaluations in the model [5]. Different hyper reduction methods are examined in the following to assess their theoretical effects on numerical stability.

### 6.1. Gappy POD and collocation methods

A simple collocation approach to hyper reduction of nonlinear mechanical systems stems directly from the Gappy POD method, which was first developed by Everson et al. [41] in the context of graphical image reconstruction from gappy data. Consider the overdetermined system of equations (44), where $\mathbf{x}$ was approximated by $\mathbf{V}\tilde{\mathbf{x}}$ (notice that we do not require an $\mathbf{M}$-orthonormal reduced basis here). Instead of restricting the residual to be orthogonal to the space spanned by the columns of $\mathbf{V}$, the collocation approach requires the residual to be exactly zero at specific collocation DoFs $i \in \mathcal{I}, |\mathcal{I}| = k$. We briefly examine two such approaches.

*a. Naïve collocation* The most simple collocation method simply requires $\mathbf{P}^T\mathbf{r}\left(\dot{\tilde{\mathbf{x}}}, \tilde{\mathbf{x}}, t\right) = \mathbf{0}$, where $\mathbf{P} \in \mathbb{R}^{m\times k}$ is a truncated permutation matrix, and $\mathbf{P}^T$ collects the $i \in \mathcal{I}$ collocation rows of $\mathbf{r}$. The collocation rows are also referred to as collocation points, or in the context of DEIM, as the DEIM "magic" points. They are typically computed using greedy algorithms [42]. The linearized system then fulfils

$$\mathbf{P}^T\mathbf{r}\left(\dot{\tilde{\mathbf{x}}}, \tilde{\mathbf{x}}, t\right) = \mathbf{0} \iff \mathbf{P}^T\left(\mathbf{M}\mathbf{V}\ddot{\tilde{\mathbf{x}}} + \mathbf{C}\mathbf{V}\dot{\tilde{\mathbf{x}}} + \mathbf{K}\mathbf{V}\tilde{\mathbf{x}} - \mathbf{f}_{\text{ext}}\right) = \mathbf{0}. \tag{77}$$

Here, we do not yet require that $\mathbf{V}$ be orthonormal to $\mathbf{M}$; the following analysis holds for any reduced basis with orthogonal columns. The sparsity properties of $\mathbf{M}$, $\mathbf{C}$ and $\mathbf{K}$ then imply that only a reduced subset $\mathcal{J}, \mathcal{I} \subset \mathcal{J}$ of the original DoFs need to be evaluated in order to compute $\mathbf{P}^T\mathbf{r}$. In particular,

- only $\mathbf{P}\mathbf{P}^T\mathbf{V}\ddot{\tilde{\mathbf{x}}}$ must be evaluated to correctly approximate the acceleration term $\mathbf{P}^T\mathbf{M}\mathbf{V}\ddot{\tilde{\mathbf{x}}}$. The full vector of accelerations $\mathbf{V}\ddot{\tilde{\mathbf{x}}}$ therefore only needs to be evaluated at the collocation rows.

- similarly, only $\mathbf{YY}^T\mathbf{V}\dot{\tilde{\mathbf{x}}}$ must be evaluated to correctly approximate the damping force term, and only $\mathbf{ZZ}^T\mathbf{V}\tilde{\mathbf{x}}$ need to be computed to approximate the forces pertaining to $\mathbf{K}$.

In a FEM context, the collocation points, as well as the rows selected by $\mathbf{Y}^T$ and $\mathbf{Z}^T$ constitute the so-called reduced mesh. It includes the nodes pertaining to the collocation points, all nodes of the adjacent elements, and any other nodes that directly influence the collocation points through contact or other phenomena. However, this approach does not generally preserve symmetry or positive (semi-) definiteness. The transformed system matrices of

$$\mathbf{P}^T\mathbf{MV}\ddot{\tilde{\mathbf{x}}} + \mathbf{P}^T\mathbf{CV}\dot{\tilde{\mathbf{x}}} + \mathbf{P}^T\mathbf{KV}\tilde{\mathbf{x}} - \mathbf{P}^T\mathbf{f}_{\text{ext}} = \mathbf{0} \tag{78}$$

$$\iff \mathbf{P}^T\mathbf{MPP}^T\mathbf{V}\ddot{\tilde{\mathbf{x}}} + \mathbf{P}^T\mathbf{CYY}^T\mathbf{V}\dot{\tilde{\mathbf{x}}} + \mathbf{P}^T\mathbf{KZZ}^T\mathbf{V}\tilde{\mathbf{x}} - \mathbf{P}^T\mathbf{f}_{\text{ext}} = \mathbf{0} \qquad \text{(sparsity)} \tag{79}$$

do not admit the application of the Poincaré separation theorem or any symmetry theorems. Hyper reduced systems based on such collocation methods can therefore generally be expected to be less stable than the FOM. We later found that this is also in line with [43], where it is shown that collocation-based methods do not preserve symmetry and therefore destroy the Lagrangian (or Hamiltonian) structure of the problem. A version of a simple collocation-based hyper reduced central difference scheme is listed below. Variants of this method have also been studied for the application in automotive crash simulation [44].

**Algorithm 2** Central difference scheme for the hyper reduced equations, using simple collocation (cf. algorithm 1)

**Input:** Previous vector of reduced displacements $\tilde{\mathbf{x}}^n$ and velocities $\dot{\tilde{\mathbf{x}}}^{n-1/2}$, nodal forces at the collocation points $\mathbf{P}^T\mathbf{f}^n_{\text{node}}$, mass matrix $\mathbf{M}^n$, time $t^n$, time step size $\Delta t$.
**Output:** Updated vector of reduced displacements $\tilde{\mathbf{x}}^{n+1}$, reduced velocities $\dot{\tilde{\mathbf{x}}}^{n+1/2}$, time $t^{n+1}$.

1: Accelerations update: $\mathbf{P}^T\ddot{\mathbf{x}}^n \leftarrow \mathbf{P}^T(\mathbf{M}^n)^{-1}\mathbf{PP}^T\mathbf{f}^n_{\text{node}}$ ▷
$\mathbf{P}^T\mathbf{f}^n_{\text{node}} = -\mathbf{P}^T\mathbf{CYY}^T\mathbf{V}\dot{\tilde{\mathbf{x}}}^{n-1/2} - \mathbf{P}^T\mathbf{KZZ}^T\mathbf{V}\tilde{\mathbf{x}}^n + \mathbf{P}^T\mathbf{f}^n_{\text{ext}}$
2: Velocity update: $\mathbf{P}^T\dot{\mathbf{x}}^{n+1/2} \leftarrow \mathbf{P}^T\dot{\mathbf{x}}^{n-1/2} + \Delta t\mathbf{P}^T\ddot{\mathbf{x}}^n$
3: Reduced displacements update: $\mathbf{P}^T\mathbf{V}\tilde{\mathbf{x}}^{n+1} \leftarrow \mathbf{P}^T\mathbf{V}\tilde{\mathbf{x}}^n + \Delta t\mathbf{P}^T\dot{\mathbf{x}}^{n+1/2}$
4: Reduced velocity update: $\dot{\tilde{\mathbf{x}}}^{n+1/2} \leftarrow \left(\tilde{\mathbf{x}}^{n+1} - \tilde{\mathbf{x}}^n\right)/\Delta t$
5: Time update: $t^n \leftarrow t^n + \Delta t$; $n \leftarrow n+1$

Alternatively, the reduced velocities can also be computed in line 2 via $\mathbf{P}^T\mathbf{V}\dot{\tilde{\mathbf{x}}}^{n+1/2} \leftarrow \mathbf{P}^T\dot{\mathbf{x}}^{n-1/2} + \Delta t\mathbf{P}^T\ddot{\mathbf{x}}^n$. Similar to the DEIM and GNAT methods, it is further possible to include more than $k$ collocation points and determine the optimum solution in the space spanned by $\mathbf{V}$ using a linear least squares solver. This typically produces more robust results [45]. It also gives rise to an online-stage error estimator, as the residual will no longer be equal to $\mathbf{0}$ for $\mathbf{P} \in \mathbb{R}^{m\times p}, p > k$.

*b. Projected collocation* The projected collocation method (cf. [46]) requires the *projection* of certain rows of the residual vector onto the reduced basis to be zero, i.e.

$$\left(\mathbf{P}^T\mathbf{V}\right)^T\mathbf{P}^T\mathbf{r}\left(\dot{\tilde{\mathbf{x}}}, \tilde{\mathbf{x}}, t\right) = \mathbf{0} \iff \mathbf{V}^T\mathbf{PP}^T\left(\mathbf{MV}\ddot{\tilde{\mathbf{x}}} + \mathbf{CV}\dot{\tilde{\mathbf{x}}} + \mathbf{KV}\tilde{\mathbf{x}} - \mathbf{f}_{\text{ext}}\right) = \mathbf{0} \tag{80}$$

if $\mathbf{V}$, with $\mathbf{V}^T\mathbf{V}$, is used as a reduced basis. In the case that $\mathbf{V}_M = \mathbf{M}^{-1/2}\mathbf{V}$ is used as a reduced basis, the equation is

$$\left(\mathbf{P}^T\mathbf{V}_M\right)^T\mathbf{P}^T\mathbf{r}\left(\dot{\tilde{\mathbf{x}}}, \tilde{\mathbf{x}}, t\right) = \mathbf{0} \iff \mathbf{V}_M^T\mathbf{PP}^T\left(\mathbf{MV}_M\ddot{\tilde{\mathbf{x}}} + \mathbf{CV}_M\dot{\tilde{\mathbf{x}}} + \mathbf{KV}_M\tilde{\mathbf{x}} - \mathbf{f}_{\text{ext}}\right) = \mathbf{0}. \tag{81}$$

As for naïve collocation, the sparsity properties of $\mathbf{M}$, $\mathbf{C}$ and $\mathbf{K}$ allow for evaluating only a reduced subset of the full-dimensional DoFs in order to correctly compute the products. Using the same notation as above yields the hyper reduced system equations

$$\mathbf{V}^T\mathbf{PP}^T\mathbf{MV}\ddot{\tilde{\mathbf{x}}} + \mathbf{V}^T\mathbf{PP}^T\mathbf{CV}\dot{\tilde{\mathbf{x}}} + \mathbf{V}^T\mathbf{PP}^T\mathbf{KV}\tilde{\mathbf{x}} - \mathbf{V}^T\mathbf{PP}^T\mathbf{f}_{\text{ext}} = \mathbf{0} \tag{82}$$

$$\iff \underbrace{\mathbf{V}^T\mathbf{PP}^T\mathbf{MPP}^T\mathbf{V}}_{\mathbf{M}_r}\ddot{\tilde{\mathbf{x}}} + \underbrace{\mathbf{V}^T\mathbf{PP}^T\mathbf{CYY}^T\mathbf{V}}_{\mathbf{C}_r}\dot{\tilde{\mathbf{x}}} + \underbrace{\mathbf{V}^T\mathbf{PP}^T\mathbf{KZZ}^T\mathbf{V}}_{\mathbf{K}_r}\tilde{\mathbf{x}} - \mathbf{V}^T\mathbf{PP}^T\mathbf{f}_{\text{ext}} = \mathbf{0}, \tag{83}$$

when $\mathbf{V}$ is used as a reduced basis. When $\mathbf{V}_M$ is used, the hyper reduced equations are described by

$$\underbrace{\mathbf{V}^T\mathbf{PP}^T\mathbf{V}}_{\mathbf{M}_r}\ddot{\tilde{\mathbf{x}}} + \underbrace{\mathbf{V}^T\mathbf{M}^{-1/2}\mathbf{PP}^T\mathbf{CYY}^T\mathbf{M}^{-1/2}\mathbf{V}}_{\mathbf{C}_r}\dot{\tilde{\mathbf{x}}} + \underbrace{\mathbf{V}^T\mathbf{M}^{-1/2}\mathbf{PP}^T\mathbf{KZZ}^T\mathbf{M}^{-1/2}\mathbf{V}}_{\mathbf{K}_r}\tilde{\mathbf{x}} - \mathbf{V}^T\mathbf{M}^{-1/2}\mathbf{PP}^T\mathbf{f}_{\text{ext}} = \mathbf{0}, \tag{84}$$

Interestingly, $\mathbf{M}_r$ inherits the symmetry and positive definiteness of $\mathbf{M}$ in either case. This is due to the fact that $\mathbf{P}\mathbf{P}^T$ is a symmetric, diagonal, idempotent matrix (corollary 4.2), since $\left(\mathbf{P}\mathbf{P}^T\right)\left(\mathbf{P}\mathbf{P}^T\right) = \mathbf{P}\mathbf{P}^T$. However, $\mathbf{C}_r$ and $\mathbf{K}_r$ do not in general inherit these properties, as $\mathbf{P} \neq \mathbf{Y}$ and $\mathbf{P} \neq \mathbf{Z}$, unless the collocation points constitute the entire mesh, in which case this method converges to the Galerkin ROM. The general stability properties and eigenvalue interlacing of the ROM are therefore not guaranteed to be preserved in this hyper reduced-order model (HROM). This is in agreement with [43], where it is shown that collocation-based methods generally do not preserve symmetry and destroy the Lagrangian structure.

We emphasize that collocation approaches can nevertheless yield satisfactory results and a high accuracy for certain problems, even though they do not preserve stability for all models. In particular, they tend to avoid over-fitting of the training simulations, which is often observed for the DEIM and GNAT methods [46, p. 98].

### 6.2. Discrete Empirical Interpolation Methods

The main difference of the Discrete Empirical Interpolation Method (DEIM) [42] and the Gauss-Newton with Approximated Tensors (GNAT) method [45, 47] compared to collocation based approaches is the approximation of the nonlinear terms using a separate reduced basis. For instance, approximating the entire nodal force vector using the DEIM yields [48]

$$\mathbf{f}_{\text{node}}(\dot{\tilde{\mathbf{x}}}, \tilde{\mathbf{x}}, t) \approx \mathbf{U}\tilde{\mathbf{f}}_{\text{node}}(\dot{\tilde{\mathbf{x}}}, \tilde{\mathbf{x}}, t), \tag{85}$$

where $\mathbf{U} \in \mathbb{R}^{m \times k_{\text{DEIM}}}$ is a reduced basis for the nodal forces, and $\tilde{\mathbf{f}}_{\text{node}}(\dot{\tilde{\mathbf{x}}}, \tilde{\mathbf{x}}, t) \in \mathbb{R}^{k_{\text{DEIM}}}$ is the vector of reduced nodal forces. The original DEIM algorithm then computes the reduced nodal forces in the online phase such that

$$\mathbf{P}^T\mathbf{U}\tilde{\mathbf{f}}_{\text{node}}(\dot{\tilde{\mathbf{x}}}, \tilde{\mathbf{x}}, t) = \mathbf{P}^T\mathbf{f}_{\text{node}}(\mathbf{V}\dot{\tilde{\mathbf{x}}}, \mathbf{V}\tilde{\mathbf{x}}, t) \iff \tilde{\mathbf{f}}_{\text{node}}(\dot{\tilde{\mathbf{x}}}, \tilde{\mathbf{x}}, t) = \left(\mathbf{P}^T\mathbf{U}\right)^{-1}\mathbf{P}^T\mathbf{f}_{\text{node}}(\mathbf{V}\dot{\tilde{\mathbf{x}}}, \mathbf{V}\tilde{\mathbf{x}}, t) \tag{86}$$

$$\mathbf{f}_{\text{node}} \approx \mathbf{U}\left(\mathbf{P}^T\mathbf{U}\right)^{-1}\mathbf{P}^T\mathbf{f}_{\text{node}}(\mathbf{V}\dot{\tilde{\mathbf{x}}}, \mathbf{V}\tilde{\mathbf{x}}, t). \tag{87}$$

Here, we do not yet require that $\mathbf{V}$ be orthonormal to $\mathbf{M}$; the following analysis holds for any reduced basis with orthogonal columns. $\mathbf{P} \in \mathbb{R}^{m \times k_{\text{DEIM}}}$ again denotes a restriction matrix. Notice that DEIM or GNAT are often only applied to the nonlinear components of the nodal forces, but they can be analogously applied to the entire residual as well. Projecting the residual once again on the reduced basis $\mathbf{V}$ yields

$$\mathbf{V}^T\left(\mathbf{M}\mathbf{V}\ddot{\tilde{\mathbf{x}}} - \mathbf{U}\left(\mathbf{P}^T\mathbf{U}\right)^{-1}\mathbf{P}^T\mathbf{f}_{\text{node}}(\mathbf{V}\dot{\tilde{\mathbf{x}}}, \mathbf{V}\tilde{\mathbf{x}}, t)\right) = \mathbf{0}. \tag{88}$$

Inserting $\mathbf{f}_{\text{node}} = \mathbf{f}_{\text{ext}} - \mathbf{C}\mathbf{V}\dot{\tilde{\mathbf{x}}} - \mathbf{K}\mathbf{V}\tilde{\mathbf{x}}$ yields the following system of hyper reduced equations

$$\mathbf{V}^T\mathbf{M}\mathbf{V}\ddot{\tilde{\mathbf{x}}} + \mathbf{V}^T\mathbf{U}\left(\mathbf{P}^T\mathbf{U}\right)^{-1}\mathbf{P}^T\mathbf{C}\mathbf{V}\dot{\tilde{\mathbf{x}}} + \mathbf{V}^T\mathbf{U}\left(\mathbf{P}^T\mathbf{U}\right)^{-1}\mathbf{P}^T\mathbf{K}\mathbf{V}\tilde{\mathbf{x}} - \mathbf{V}^T\mathbf{U}\left(\mathbf{P}^T\mathbf{U}\right)^{-1}\mathbf{P}^T\mathbf{f}_{\text{ext}} = \mathbf{0}. \tag{89}$$

Introducing the same reduced mesh notation with the restriction matrices $\mathbf{Y}$ and $\mathbf{Z}$ as in equation (78) gives the following final representation

$$\mathbf{V}^T\mathbf{M}\mathbf{V}\ddot{\tilde{\mathbf{x}}} + \mathbf{V}^T\mathbf{U}\left(\mathbf{P}^T\mathbf{U}\right)^{-1}\mathbf{P}^T\mathbf{C}\mathbf{Y}\mathbf{Y}^T\mathbf{V}\dot{\tilde{\mathbf{x}}} + \mathbf{V}^T\mathbf{U}\left(\mathbf{P}^T\mathbf{U}\right)^{-1}\mathbf{P}^T\mathbf{K}\mathbf{Z}\mathbf{Z}^T\mathbf{V}\tilde{\mathbf{x}} - \mathbf{V}^T\mathbf{U}\left(\mathbf{P}^T\mathbf{U}\right)^{-1}\mathbf{P}^T\mathbf{f}_{\text{ext}} = \mathbf{0}, \tag{90}$$

or equivalently if $\mathbf{V}_M$ is used as a reduced basis:

$$\begin{aligned}\ddot{\tilde{\mathbf{x}}} &+ \mathbf{V}^T\mathbf{M}^{-1/2}\mathbf{U}\left(\mathbf{P}^T\mathbf{U}\right)^{-1}\mathbf{P}^T\mathbf{C}\mathbf{Y}\mathbf{Y}^T\mathbf{M}^{-1/2}\mathbf{V}\dot{\tilde{\mathbf{x}}} + \\ &+ \mathbf{V}^T\mathbf{M}^{-1/2}\mathbf{U}\left(\mathbf{P}^T\mathbf{U}\right)^{-1}\mathbf{P}^T\mathbf{K}\mathbf{Z}\mathbf{Z}^T\mathbf{M}^{-1/2}\mathbf{V}\tilde{\mathbf{x}} - \mathbf{V}^T\mathbf{M}^{-1/2}\mathbf{U}\left(\mathbf{P}^T\mathbf{U}\right)^{-1}\mathbf{P}^T\mathbf{f}_{\text{ext}} = \mathbf{0}.\end{aligned} \tag{91}$$

The GNAT method can be seen as a generalization of DEIM, in which the restriction matrix $\mathbf{P}^T \in \mathbb{R}^{m \times p}$ may select more than $k_{\text{DEIM}}$ rows (i.e. $p \geq k_{\text{DEIM}}$) [46]. The equations are then solved in a least squares sense, typically yielding more robust approximations. The GNAT-generalized form of (90) is

$$\mathbf{V}^T\mathbf{M}\mathbf{V}\ddot{\tilde{\mathbf{x}}} + \mathbf{V}^T\mathbf{U}\left(\mathbf{P}^T\mathbf{U}\right)^{\dagger}\mathbf{P}^T\mathbf{C}\mathbf{Y}\mathbf{Y}^T\mathbf{V}\dot{\tilde{\mathbf{x}}} + \mathbf{V}^T\mathbf{U}\left(\mathbf{P}^T\mathbf{U}\right)^{\dagger}\mathbf{P}^T\mathbf{K}\mathbf{Z}\mathbf{Z}^T\mathbf{V}\tilde{\mathbf{x}} - \mathbf{V}^T\mathbf{U}\left(\mathbf{P}^T\mathbf{U}\right)^{\dagger}\mathbf{P}^T\mathbf{f}_{\text{ext}} = \mathbf{0}, \tag{92}$$

where $(\cdot)^{\dagger}$ denotes the Moore-Penrose pseudoinverse. Equation (91) can be written analogously when an $\mathbf{M}$-orthonormal reduced basis is used for GNAT.

As for the collocation methods, no general statement about the symmetry or eigenvalues of the hyper reduced system matrices can be inferred. $\mathbf{V}^T\mathbf{M}\mathbf{V}$ preserves symmetry and positive definiteness; however the hyper reduced damping and stiffness matrices do not enjoy general symmetry or eigenvalue interlacing properties. This indicates that the DEIM method may be prone to inherent instability, and cannot generally be expected to have larger stability time steps than the FOM. In fact, several authors have reported instabilities of DEIM-type hyper-reduced models, or observed that DEIM may produce nonsymmetric stiffness matrices [5, 49–52]. These instabilities are problem-dependent and also partly owing to the fact that the focus in the conception of the DEIM and GNAT methods was primarily on accuracy and not numerical stability [5].

As for collocation methods, we later found that it is shown in [43] that DEIM and GNAT destroy symmetry and the Lagrangian structure of the problem. Carlberg et al. also describe further approaches to preserve Lagrangian structure in the presence of hyper reduction using reduced basis sparsification and matrix gappy POD, with applications to structural dynamics [43].

### 6.3. ECSW hyper reduction method

The ECSW method was developed for the reduction of second-order nonlinear dynamical systems in the context of structural dynamics FEM models [5, 53]. The concept is related to the quadrature method for computer graphics initially described by [54]. It is different to the aforementioned hyper reduction methods in that it does not directly approximate any nonlinear force terms using a gappy completion approach. Instead, it approximates directly the associated *projected* force vectors, which can be interpreted physically as an energy. As DEIM, it is typically only applied to the nonlinear terms which are unaffordable to compute on the FOM dimension. For the stability analysis, consider the application of ECSW to the entire vector of nodal forces $\mathbf{f}_{\text{node}} \in \mathbb{R}^m$ with an $\mathbf{M}$-orthonormal reduced basis $\mathbf{V}_M = \mathbf{M}^{-1/2}\mathbf{V}$

$$\mathbf{V}_M^T \mathbf{M} \mathbf{V}_M \ddot{\tilde{\mathbf{x}}} - \mathbf{V}_M^T \mathbf{f}_{\text{node}} \left( \mathbf{V}_M \dot{\tilde{\mathbf{x}}}, \mathbf{V}_M \tilde{\mathbf{x}}, t \right) = 0 \tag{93}$$

$$\ddot{\tilde{\mathbf{x}}} - \mathbf{V}_M^T \mathbf{f}_{\text{node}} \left( \mathbf{V}_M \dot{\tilde{\mathbf{x}}}, \mathbf{V}_M \tilde{\mathbf{x}}, t \right) = 0. \tag{94}$$

The ECSW method reduces the number of nonlinear function evaluations by introducing element weighting factors in the nodal force assembly process at the element level. It approximates the projection of the nonlinear forces onto $\mathbf{V}$ by

$$\mathbf{V}_M^T \mathbf{f}_{\text{node}} \left( \mathbf{V}_M \dot{\tilde{\mathbf{x}}}, \mathbf{V}_M \tilde{\mathbf{x}}, t \right) = \mathbf{V}_M^T \sum_{e \in \Omega} \mathbf{L}_e^T \mathbf{f}_e \left( \mathbf{L}_e \mathbf{V}_M \dot{\tilde{\mathbf{x}}}, \mathbf{L}_e \mathbf{V}_M \tilde{\mathbf{x}}, t \right) \approx \mathbf{V}_M^T \sum_{e \in \Omega} \xi_e^* \mathbf{L}_e^T \mathbf{f}_e \left( \mathbf{L}_e \mathbf{V}_M \dot{\tilde{\mathbf{x}}}, \mathbf{L}_e \mathbf{V}_M \tilde{\mathbf{x}}, t \right). \tag{95}$$

Here, the sum denotes a summation over all elements $e$ within the FE mesh $\Omega$, and $\mathbf{L}_e \in \mathbb{R}^{m_e \times m}$ is the element connectivity matrix which ensures correct assembly of the nodal forces for each index (cf. 3.4). $m_e$ is the number of DoF in the element $e$, $\xi_e \geq 0$ is the non-negative weighting factor introduced by the ECSW method. The goal of the ECSW training stage is to select the $\xi_e^*$ such that a high approximation accuracy is achieved *and* as many of the $\xi_e$ as possible will be zero. This is equivalent to a minimization problem using the zero norm [5]. The benefit from this approximation is that the element contributions with zero weighting factors can be elegantly omitted from the sum, giving rise to the approximation

$$\mathbf{V}_M^T \mathbf{f}_{\text{node}} \left( \mathbf{V}_M \dot{\tilde{\mathbf{x}}}, \mathbf{V}_M \tilde{\mathbf{x}}, t \right) \approx \mathbf{V}_M^T \sum_{e \in \Omega_{\text{RM}}} \xi_e^* \mathbf{L}_e^T \mathbf{f}_e \left( \mathbf{L}_e \mathbf{V}_M \dot{\tilde{\mathbf{x}}}, \mathbf{L}_e \mathbf{V}_M \tilde{\mathbf{x}}, t \right), \tag{96}$$

where $\Omega_{\text{RM}}$ is the reduced mesh which only includes the elements with $\xi_e > 0$. Exploiting $\mathbf{f}_{\text{node}} = \mathbf{f}_{\text{ext}} - \mathbf{C}\mathbf{V}_M \dot{\tilde{\mathbf{x}}} - \mathbf{K}\mathbf{V}_M \tilde{\mathbf{x}}$, the stiffness and damping matrices can be reduced in a similar manner. Using ECSW, the linearized hyper reduced equations hence become

$$\underbrace{\mathbf{I}}_{\mathbf{M}_r} \ddot{\tilde{\mathbf{x}}} + \underbrace{\mathbf{V}_M^T \left( \sum_{e \in \Omega_{\text{RM}}} \xi_e^* \mathbf{L}_e^T \mathbf{C}_e \mathbf{L}_e \right) \mathbf{V}_M}_{\mathbf{C}_r} \dot{\tilde{\mathbf{x}}} + \underbrace{\mathbf{V}_M^T \left( \sum_{e \in \Omega_{\text{RM}}} \xi_e^* \mathbf{L}_e^T \mathbf{K}_e \mathbf{L}_e \right) \mathbf{V}_M}_{\mathbf{K}_r} \tilde{\mathbf{x}} - \mathbf{V}_M^T \mathbf{f}_{\text{ext}} = 0. \tag{97}$$

$\mathbf{M}_r = \mathbf{I} = \mathbf{V}_M^T \mathbf{M} \mathbf{V}_M$ naturally preserves symmetry and positive definiteness of $\mathbf{M}$. Further, the terms

$$\sum_{e \in \Omega_{\text{RM}}} \xi_e^* \mathbf{L}_e^T \mathbf{C}_e \mathbf{L}_e \qquad \text{and} \qquad \sum_{e \in \Omega_{\text{RM}}} \xi_e^* \mathbf{L}_e^T \mathbf{K}_e \mathbf{L}_e$$

constitute sums of symmetric and positive semi-definite matrices, with positive weighting factors $\xi_e^*$. The sums are therefore also symmetric and positive semi-definite. By the analysis in section 5, both $\mathbf{C}_r$ and $\mathbf{K}_r$ therefore preserve the symmetry and positive semi-definiteness of $\mathbf{C}$ and $\mathbf{K}$. These are important properties for the stability and physical meaningfulness of the ECSW hyper reduction method. By the same analysis as in section 5, the stability time step for the central difference time integration scheme and the ECSW hyper reduced model is determined by (theorem 5.5)

$$\Delta t_{\mathrm{crit},k} = \frac{2}{\sqrt{\tilde{\mu}_k^{HROM}}}\left(\sqrt{\tilde{\xi}_k^2+1}-\tilde{\xi}_k\right), \qquad \tilde{\xi}_k = \frac{a_1}{2\sqrt{\tilde{\mu}_k^{HROM}}} + \frac{a_2\sqrt{\tilde{\mu}_k^{HROM}}}{2}, \tag{98}$$

where $\tilde{\mu}_k^{HROM}$ is the largest eigenvalue of the matrix $\mathbf{M}_r^{-1}\mathbf{K}_r$. Further, $\mathbf{M}_r^{-1}\mathbf{K}_r$ is given by

$$\mathbf{M}_r^{-1}\mathbf{K}_r = \mathbf{V}^T\mathbf{M}^{-1/2}\left(\sum_{e\in\Omega_{\mathrm{RM}}} \xi_e^*\mathbf{L}_e^T\mathbf{K}_e\mathbf{L}_e\right)\mathbf{M}^{-1/2}\mathbf{V}. \tag{99}$$

The analysis of the eigenvalues of $\mathbf{M}_r^{-1}\mathbf{K}_r$ for ECSW hyper reduction is somewhat complicated by the contrary influences of two separate phenomena (figure 2). As is known from Poincaré's eigenvalue separation theorem and section 5, the two-sided multiplication of $\mathbf{V}^T(\cdot)\mathbf{V}$ tends to reduce the magnitude of the largest eigenvalue. The projection onto the – typically – lower-frequent modes obtained from an SVD of the snapshot matrix tends to reduce the largest eigenvalue even further (cf. section 7 for a short analysis of the influence of the modal content of $\mathbf{V}$). However, the approximation of $\mathbf{K} \approx \sum_{e\in\Omega_{\mathrm{RM}}} \xi_e^*\mathbf{L}_e^T\mathbf{K}_e\mathbf{L}_e$ can potentially increase the magnitude of the largest eigenvalue for $\xi_e^* > 1$. By theorem 3.4, the largest eigenvalue in absolute magnitude of $\mathbf{M}_r^{-1}\mathbf{K}_r$ can be bounded by

$$|\lambda_{\max}^{\mathrm{HROM}}| \leq \max_e |\lambda_{\max}(\mathbf{M}_e^{-1}\xi_e^*\mathbf{K}_e)| = \max_e |\xi_e^*\lambda_e^{\mathrm{FOM}}|. \tag{100}$$

Furthermore, it is to be expected that many $\xi_e^* > 1$, in order to compensate for the contribution of the omitted elements.

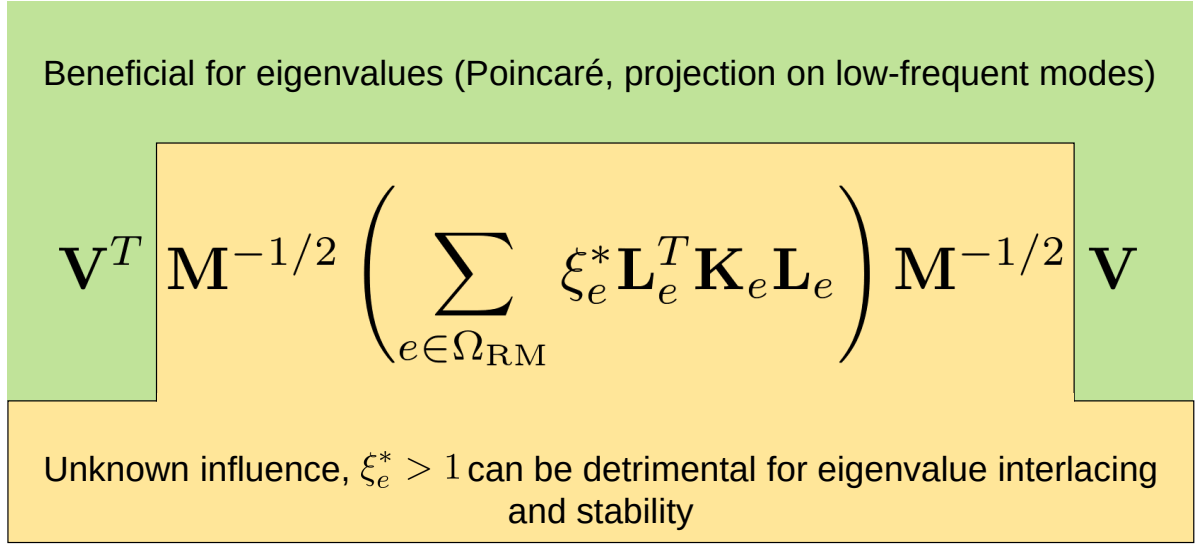

FIG. 2: Opposing influences on the eigenvalues of the ECSW hyper reduced system matrix $\mathbf{M}_r^{-1}\mathbf{K}_r$.

## 7. INFLUENCE OF THE MODAL CONTENT OF THE REDUCED BASIS

The analysis so far has not made any assumptions on the actual content of the reduced basis and is therefore generally valid regardless of which modes are selected or where the truncation occurs. However, it turns out that the content of $\mathbf{V}$ or $\mathbf{V}_M$ is important for significantly increasing the critical time step beyond that of the FOM. If the basis vectors are selected appropriately, they will be (approximately) orthogonal to the modes with the highest frequency and hence the largest eigenvalues, such that their influence will be removed by projecting onto $\mathbf{V}$. A suitable reduced basis therefore filters out the highest frequencies and largest eigenvalues of $\mathbf{M}^{-1}\mathbf{K}$. Consider the simple example of an oscillating 1D string, e.g. a guitar string, with uniform properties, discretized into $m-1$ elements and $m$ nodes (figure 3). Its mass matrix and stiffness matrices are then given by

$$\mathbf{M} = M \begin{bmatrix} 1/2 & & & & & \\ & 1 & & & & \\ & & \ddots & & & \\ & & & \ddots & & \\ & & & & 1 & \\ & & & & & 1/2 \end{bmatrix} \in \mathbb{R}^{m\times m}; \quad \mathbf{K} = K \begin{bmatrix} 100 & -1 & & & & \\ -1 & 2 & -1 & & & \\ & & \ddots & & & \\ & & & \ddots & & \\ & & & -1 & 2 & -1 \\ & & & & -1 & 100 \end{bmatrix} \in \mathbb{R}^{m\times m}. \tag{101}$$

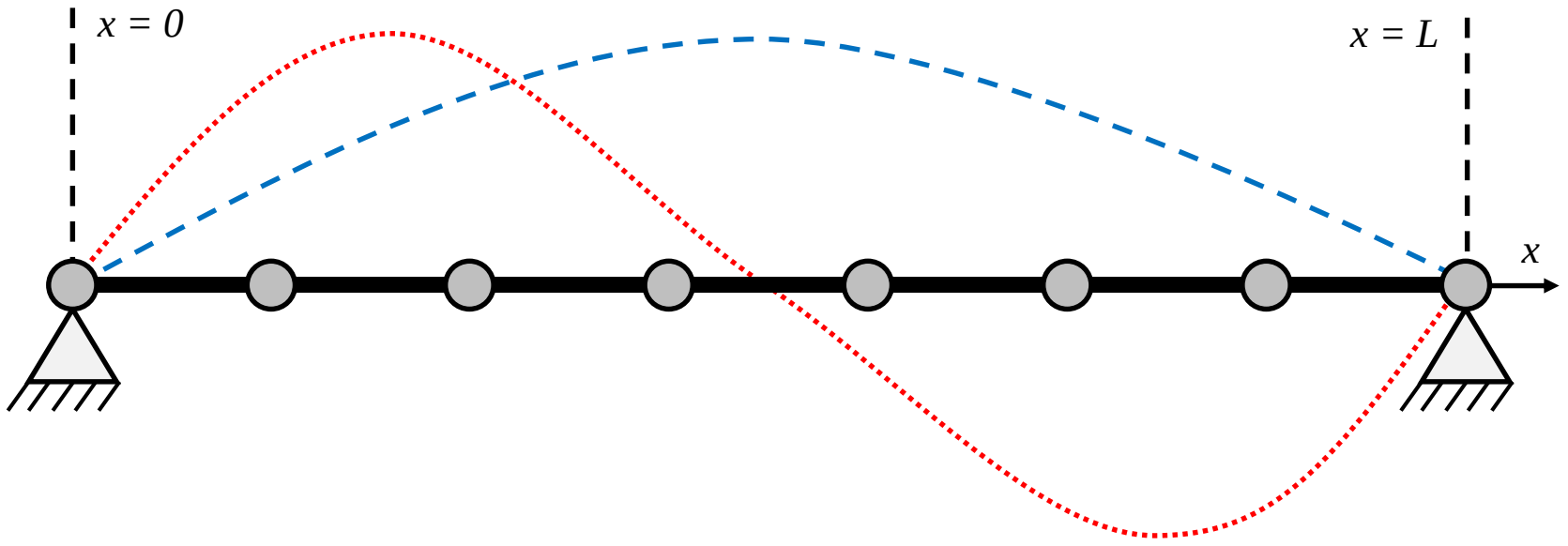

FIG. 3: Guitar string model for a discretization using $m = 8$ nodes.

Here, the boundary conditions are enforced by adding stiffness at $K_{1,1}$ and $K_{m,m}$. The matrix $\mathbf{M}^{-1}\mathbf{K}$ is then

$$\mathbf{M}^{-1}\mathbf{K} = \frac{K}{M} \begin{bmatrix} 2 & & & & & \\ & 1 & & & & \\ & & \ddots & & & \\ & & & \ddots & & \\ & & & & 1 & \\ & & & & & 2 \end{bmatrix} \begin{bmatrix} 100 & -1 & & & & \\ -1 & 2 & -1 & & & \\ & & \ddots & & & \\ & & & \ddots & & \\ & & & -1 & 2 & -1 \\ & & & & -1 & 100 \end{bmatrix} = \frac{K}{M} \begin{bmatrix} 200 & -2 & & & & \\ -1 & 2 & -1 & & & \\ & & \ddots & & & \\ & & & \ddots & & \\ & & & -1 & 2 & -1 \\ & & & & -2 & 200 \end{bmatrix} \tag{102}$$

For a discretization with $m = 100$ nodes, a length $L = 1$, $M = 1$ and $K = 10$, figure 4 visualizes the first three and the last three modes (eigenvectors) of $\mathbf{M}^{-1}\mathbf{K}$, along with their corresponding eigenvalues. Notice that the most "noisy" mode is the one with the highest frequency and the highest eigenvalue, and it also determines the stability time step for explicit analysis.
Now consider a Galerkin ROM where the reduced basis has captured the dominant system behaviour. Normally this is done from an SVD of collected snapshots; however let us first suppose for the sake of simplicity that the string was excited at relatively low frequencies and that $\mathbf{V}$ therefore contains the 10 first eigenvectors of $\mathbf{M}^{-1}\mathbf{K}$. The Galerkin ROM system matrix $\mathbf{M}_r^{-1}\mathbf{K}_r$ then exhibits natural eigenvalue interlacing as discussed in section 5.[5] In addition, the high-frequent modes are mostly orthogonal to the modes in the reduced basis, which is why they are filtered out. Figure 5 visualizes the first three and the last three modes of the Galerkin ROM. Notice that the largest eigenvalue of the ROM is in this case about 2000 times smaller than the largest eigenvalue of the FOM. In the absence of damping, this means that the critical time step of the central difference scheme is $\left(1/\sqrt{2000}\right)^{-1} \approx 44.72$ times larger (cf. section 5) than that of the FOM.

Next, let us suppose that the system is mainly excited at the very high frequencies, and the snapshots can be well approximated using a linear combination of the 10 last modes of the FOM. If $\mathbf{V}$ contains the 10 last modes of $\mathbf{M}^{-1}\mathbf{K}$, eigenvalue interlacing can still be expected; however the modes with the highest eigenvalues will no longer be orthogonal to $\mathbf{V}$. It is therefore to be expected that the largest eigenvalues of the Galerkin ROM are comparable to those of the FOM, and that the critical time step size does not change much. Indeed, this is also observed in practice (figure 6). A significant increase of the critical time step is therefore strongly dependent on the modal content of $\mathbf{V}$, and hinges on whether the high-frequent contributions are captured by the reduced basis $\mathbf{V}$.

[5] One must be careful not to modify the space spanned by the basis vectors when performing the $\mathbf{M}$-orthonormalization. We follow the approach of [55]. Notice that for every $\mathbf{M}$-orthonormal basis $\mathbf{V}_M$ and diagonal positive definite $\mathbf{M}$, there is an orthonormal basis $\mathbf{V} = \mathbf{M}^{1/2}\mathbf{V}_M$, as $\mathbf{I} = \mathbf{V}_M^T\mathbf{M}\mathbf{V}_M = \mathbf{V}^T\mathbf{M}^{-1/2}\mathbf{M}\mathbf{M}^{-1/2}\mathbf{V} = \mathbf{V}^T\mathbf{V}$.

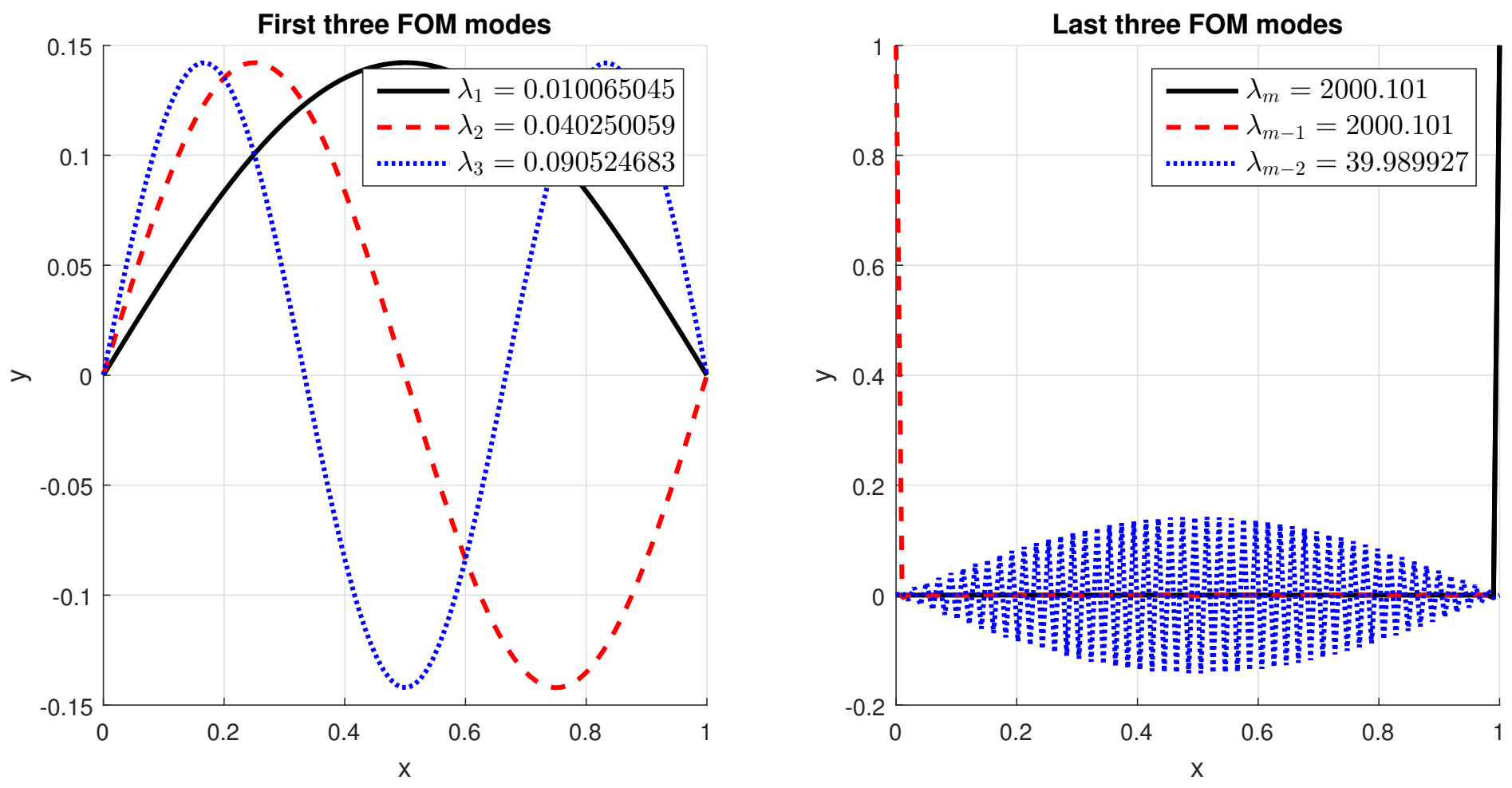

FIG. 4: Oscillation modes of the full-order guitar string model, with $m = 100$, $M = 1$, $K = 10$, and $L = 1$. The corresponding eigenvalues of $\mathbf{M}^{-1}\mathbf{K}$ are denoted as $\lambda_i$.

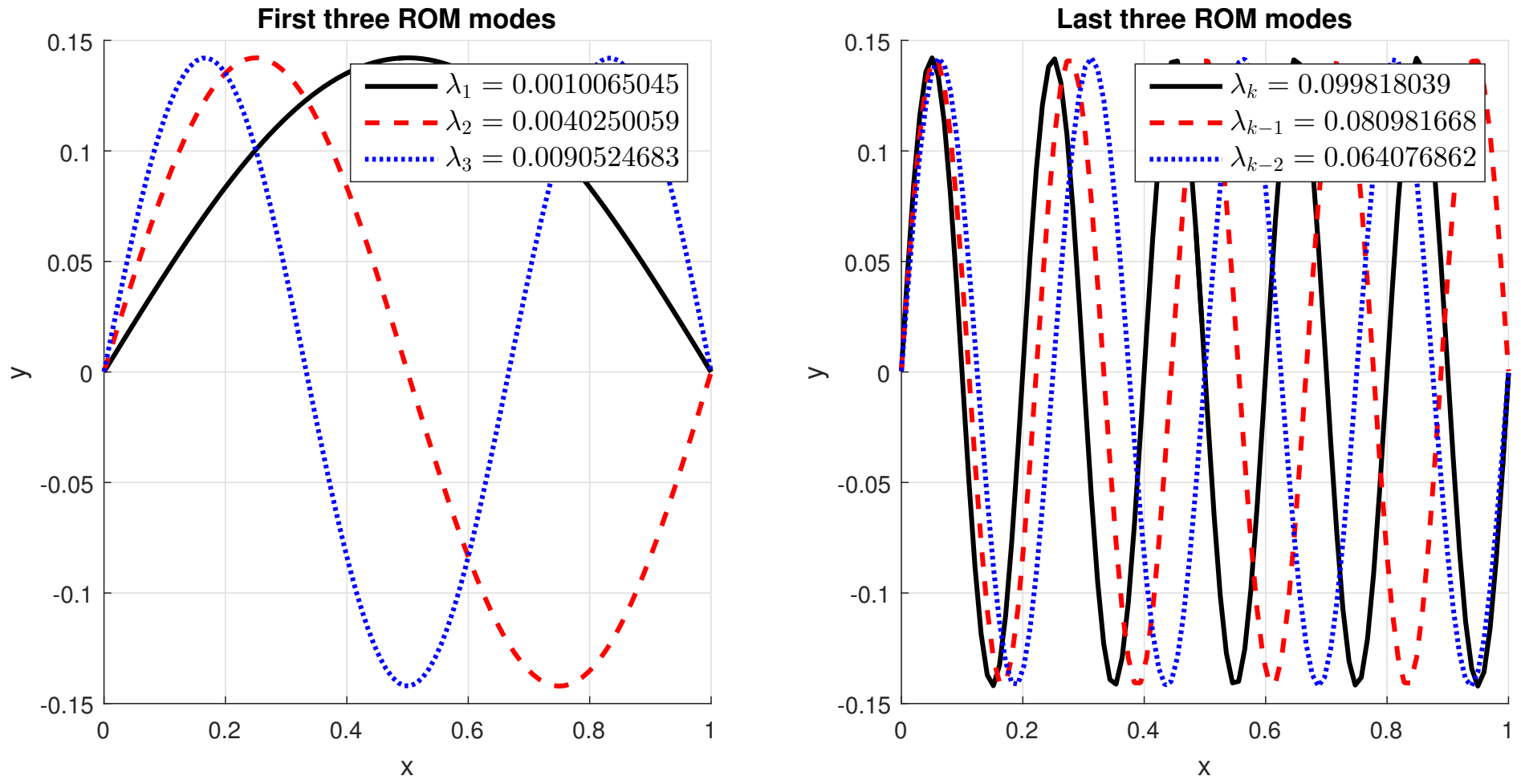

FIG. 5: Oscillation modes of the Galerkin ROM guitar string model, with $m = 100$, $M = 1$, $K = 10$, and $L = 1$. Here, the reduced basis contains the $k = 10$ first eigenvectors of $\mathbf{M}^{-1}\mathbf{K}$. The corresponding eigenvalues of $\mathbf{M}_r^{-1}\mathbf{K}_r$ are denoted as $\lambda_i$.

Finally, let us briefly analyze the impact of ECSW hyper-reduction on the stability time step using a simple example with only $m = 5$ nodes. Let us further suppose that only the second oscillation mode was observed during training simulations, and the ECSW algorithm has therefore sampled the second element (figure 7, shown in green).

Due to symmetry considerations, we further assume that the weighting factor of the sampled element, as determined by the ECSW algorithm, is $\xi_2^* = 4$. The reader can verify that, using the same notation as in (101), the FOM system matrix $\mathbf{M}^{-1}\mathbf{K}$ is given by

$$\mathbf{M}^{-1}\mathbf{K} = \frac{K}{M}\begin{bmatrix} 200 & -2 & & & \\ -1 & 2 & -1 & & \\ & -1 & 2 & -1 & \\ & & -1 & 2 & -1 \\ & & & -2 & -200 \end{bmatrix}. \tag{103}$$

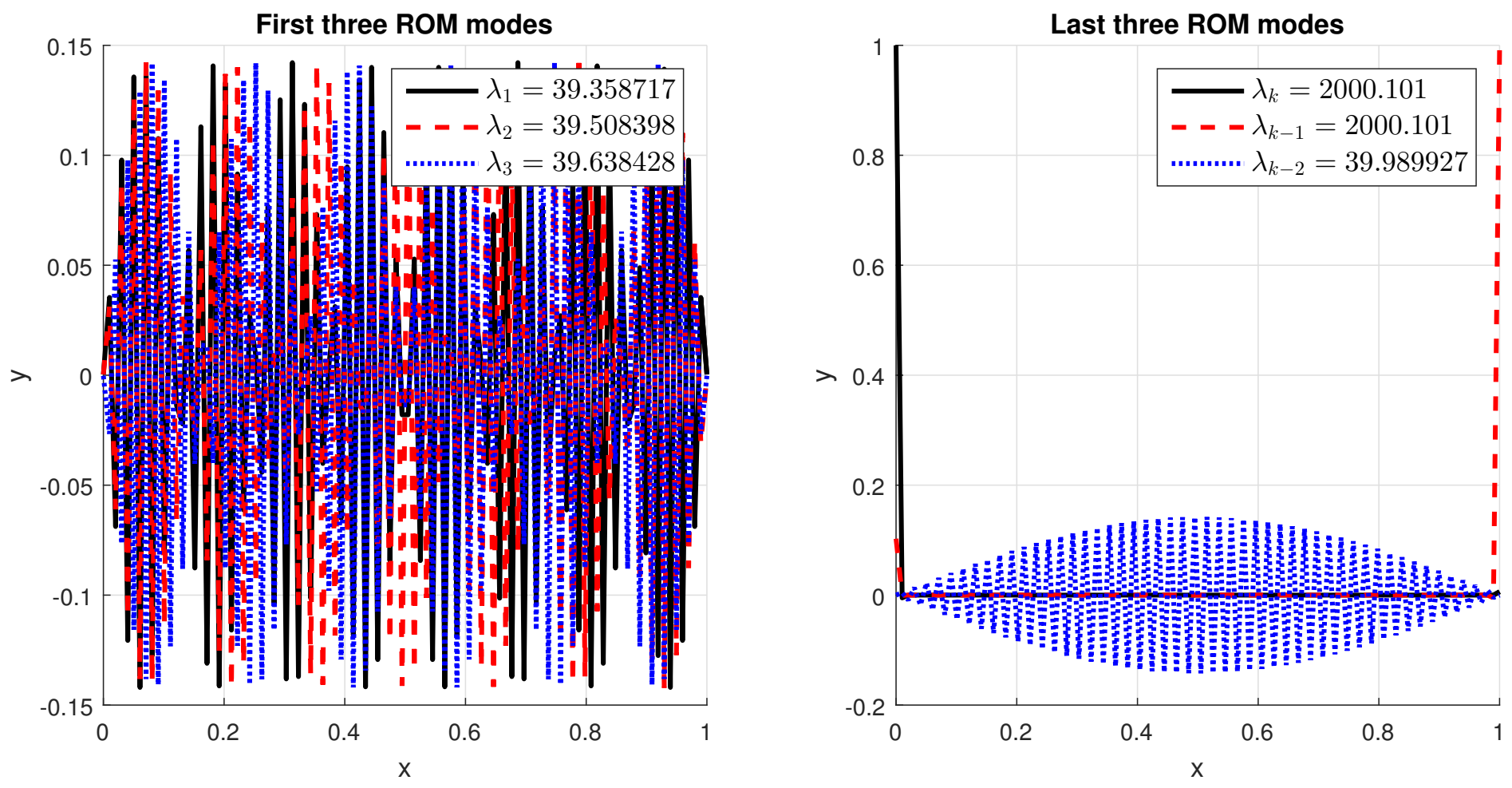

FIG. 6: Oscillation modes of the Galerkin ROM guitar string model, with $m = 100$, $M = 1$, $K = 10$, and $L = 1$. Here the reduced basis contains the $k = 10$ last eigenvectors of $\mathbf{M}^{-1}\mathbf{K}$. The corresponding eigenvalues of $\mathbf{M}_r^{-1}\mathbf{K}_r$ are denoted as $\lambda_i$.

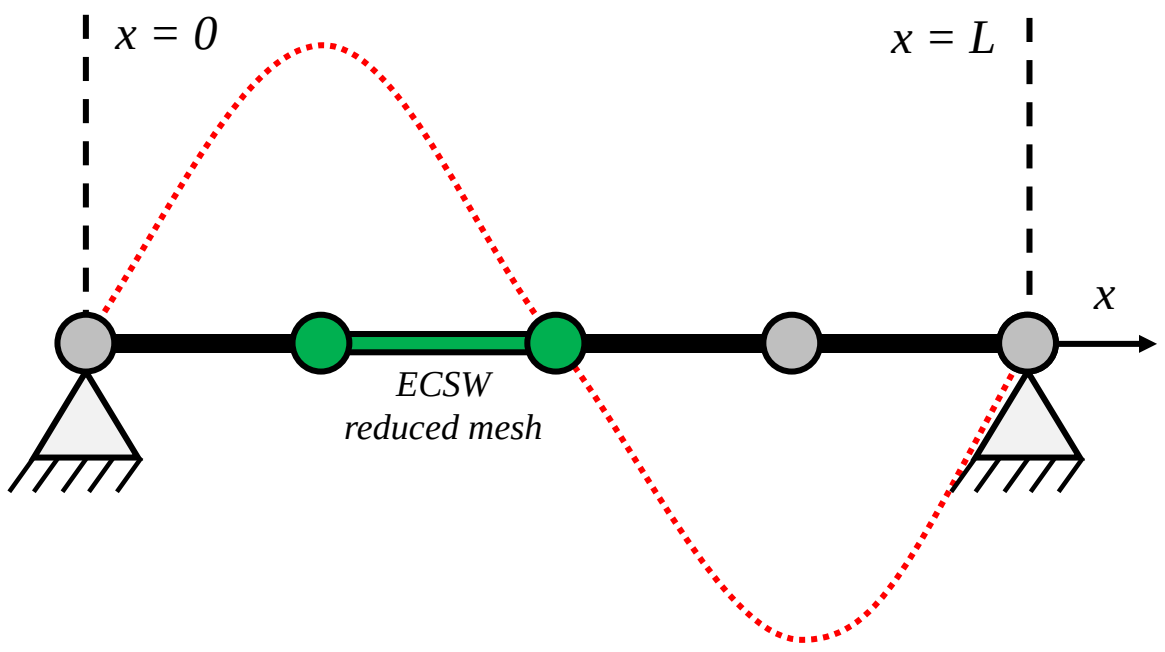

FIG. 7: ECSW hyper reduced guitar string model for $m = 5$ nodes. The reduced mesh is highlighted in green.

Again assuming $M = 1$ and $K = 10$, the eigenvalues of the FOM are

$$\left[\lambda_1^{FOM}, \lambda_2^{FOM}, \lambda_3^{FOM}, \lambda_4^{FOM}, \lambda_5^{FOM}\right] = [5.81, 19.90, 34.09, 2000.101, 2000.101]\,. \tag{104}$$

The hyper reduced system matrix of this example problem is then given by (99)

$$\mathbf{M}_r^{-1}\mathbf{K}_r = \mathbf{V}^T\mathbf{M}^{-1/2}\left(\xi_2^*\mathbf{L}_2^T\mathbf{K}_2\mathbf{L}_2\right)\mathbf{M}^{-1/2}\mathbf{V} = \mathbf{V}^T \underbrace{\left(\frac{4K}{M}\begin{bmatrix} 0 & 0 & 0 & 0 & 0 \\ 0 & 1 & -1 & 0 & 0 \\ 0 & -1 & 1 & 0 & 0 \\ 0 & 0 & 0 & 0 & 0 \\ 0 & 0 & 0 & 0 & 0 \end{bmatrix}\right)}_{=\mathbf{R}} \mathbf{V} \tag{105}$$

Here, $\mathbf{R}$ is a symmetric and positive semi-definite matrix with eigenvalues

$$\left[\lambda_1^R, \lambda_2^R, \lambda_3^R, \lambda_4^R, \lambda_5^R\right] = [0, 0, 0, 0, 80]\,. \tag{106}$$

Applying $\mathbf{V}^T(\cdot)\,\mathbf{V}$, where in this case we suppose that $\mathbf{V}$ contains only the second mode, yields the eigenvalues of the ECSW hyper reduced system. As $\mathbf{V}$ contains only one basis vector, the only eigenvalue is

$$\lambda^{HROM} = 20. \tag{107}$$

Conversely, a simple Galerkin projection without hyper reduction would have yielded a largest eigenvalue of $\lambda^{ROM} = 19.9$. The ECSW hyper-reduction scheme benefits from larger stable time step sizes than the FOM if $\mathbf{V}$ removes high-frequent contributions, and/or if critical elements are removed by the sampling method so that their contributions need not be evaluated anymore during the online stage. This can be observed in the above example, where the two critical elements at the left and the right boundaries are not part of the reduced mesh any longer. If, however, the critical elements had been included in the reduced mesh (e.g. by sampling element 1 instead of element 2), their contribution would have been amplified by the weighting factor, with an adverse impact on the stable time step unless these high-frequent modes are removed by the projection onto $\mathbf{V}$.

## 8. CONCLUSIONS

Conditions for the numerical stability of explicit integration methods for Galerkin ROMs and different hyper reduction methods have been derived. It was proved, for the first time as far as the authors are aware, that Galerkin ROMs with an $\mathbf{M}$-orthonormal reduced basis guarantee that the ROM critical time step cannot drop below that of the FOM. Further, a formula was derived which enables the computation of the ROM critical time step in the online simulation of Galerkin ROMs integrated using the central difference method. It was shown that the modal content of the reduced basis also plays an important role in improving the stability, and that the choice of low-frequent basis vectors can significantly improve the numerical stability, as these vectors will likely be (approximately) orthogonal to high-frequent modes with large eigenvalues. Both findings strongly support the understanding of why Galerkin ROMs have often been found to exhibit larger stability time steps. Further, various hyper reduction methods were analyzed for numerical stability. It was found that all analyzed node-based hyper reduction approaches, such as collocation, DEIM, and GNAT, do not generally preserve important symmetry or positive (semi-) definiteness properties of the involved system matrices for second-order nonlinear dynamical systems. It was also found that, as described first by [5, 53], the ECSW hyper reduction method preserves symmetry and positive (semi-) definiteness of the system matrices, which greatly benefits numerical stability. Two opposing mechanisms have been found to impact the critical stability time step of explicitly integrated ECSW hyper reduced models in this work. The first mechanism is the Galerkin projection $\mathbf{V}^T(\cdot)\mathbf{V}$, which tends to filter out the high-frequent modes through projection, and guarantees eigenvalue interlacing. The second, contrary mechanism is due to the introduction of force weighting factors $\xi_e^*$ at the element level, which can potentially reduce the critical stability time step compared to the FOM for $\xi_e^* > 1$. A formula was presented, for the first time as far as the authors are aware, for computing the critical stability time step of ECSW hyper reduced models equipped with a central difference time integration scheme. There is reason to hypothesize that the beneficial contribution of the Galerkin projection often dominates in practice, therefore enabling larger time steps even for hyper reduced models (cf. [5]).

We emphasize that the main focus of this work was numerical stability. Consistency and convergence to the true solution are two other important aspects which need to be satisfied in order to obtain sufficiently accurate results. Both projection-based ROMs and hyper reduction methods introduce new consistency errors in addition to those that may already be present in the original FEM model. It is currently "hoped" that the consistency errors remain small enough if the training simulations are well-chosen. For linear problems, the Lax equivalence theorem (Lax-Richtmyer-theorem) [56] typically ensures that, if a numerical method is stable and the approximation is consistent, then the method is also convergent. Further, we have applied linear stability theory to a nonlinear problem, around a linearization point. This is in line with the equivalent analysis for classical explicit FEM [9], but the linearization may not always remain valid for large time step sizes. This is often observed, for example, in case of contact problems. There are therefore additional problem dependent limits on the maximum admissible time step.

[1] E. Gutiérrez , J. J. L. Cela. Improving explicit time integration by modal truncation techniques. *Earthquake Engineering & Structural Dynamics.* 1998; 27(12):1541–1557.

[2] P. Krysl, S. Lall, J. E. Marsden. Dimensional model reduction in non-linear finite element dynamics of solids and structures. *International Journal for Numerical Methods in Engineering.* 2001; 51:479–504.

[3] C. Bucher. Stabilization of explicit time integration by modal reduction. In: W. A. Wall, K.-U. Bletzinger, K. Schweizerhof, eds. *Trends in Computational Structural Mechanics*, CIMNE; 2001; Barcelona, Spain.

[4] Z. A. Taylor, S. Crozier, S. Ourselin. A Reduced Order Explicit Dynamic Finite Element Algorithm for Surgical Simulation. *IEEE Transactions on Medical Imaging.* 2011; 30(9):1713–1721.

[5] C. Farhat, P. Avery, T. Chapman, J. Cortial. Dimensional reduction of nonlinear finite element dynamic models with finite rotations and energy-based mesh sampling and weighting for computational efficiency. *International Journal for Numerical Methods in Engineering.* 2014; 98(9):625–662.

[6] F. Bamer, A. K. Amiri, C. Bucher. A new model order reduction strategy adapted to nonlinear problems in earthquake engineering. *Earthquake Engineering & Structural Dynamics.* 2017; 46:537–559.

[7] N. M. Newmark. A Method of Computation for Structural Dynamics. *Journal of the Engineering Mechanics Division.* 1959; EM 3:67–94.
[8] K. J. Bathe. *Finite Element Procedures.* 2nd ed. 2014.
[9] T. Belytschko, W. K. Liu, B. Moran, K. I. Elkhodary. *Nonlinear Finite Elements for Continua and Structures, Second Edition.* Wiley; Chichester, UK; 2014.
[10] Livermore Software Technology Corporation (LSTC). *LS-DYNA R8.0 Keyword User's Manual.* 2015.
[11] Dassault Systèmes Simulia. *Abaqus Keywords Reference Manual (6.12).* Dassault Systèmes; 2012.
[12] Livermore Software Technology Corporation (LSTC). *LS-DYNA Theory Manual.* 2006.
[13] A. M. Lyapunov. Problème général de la stabilité du mouvement. *Annales de la faculté des sciences de Toulouse.* 1907; 2(9):203–474.
[14] D. G. Dahlquist. A special stability problem for linear multistep methods. *BIT Numerical Mathematics.* 1963; 3(1):27–43.
[15] E. J. Routh. *A treatise on the stability of a given state of motion, particularly steady motion.* London, British Empire: Macmillan and Co.; 1877.
[16] A. Hurwitz. Ueber die Bedingungen, unter welchen eine Gleichung nur Wurzeln mit negativen reellen Theilen besitzt. *Mathematische Annalen.* 1895; 46(2):273–284.
[17] N. Halko, P.-G. Martinsson, J. A. Tropp. Finding Structure with Randomness: Probabilistic Algorithms for Constructing Approximate Matrix Decompositions. *SIAM Review.* 2011; 53(2):217–288.
[18] J. R. Magnus, H. Neudecker. *Matrix Differential Calculus with Applications in Statistics and Econometrics.* Wiley; Chichester, UK; 3rd ed. 2007.
[19] T. Belytschko, P. Smolinski, W. K. Liu. Stability of multi-time step partitioned integrators for first-order finite element systems. *Computer Methods in Applied Mechanics and Engineering.* 1985; 49(3):281–297.
[20] I. Fried. Discretization and Round-Off Errors in the Finite Element Analysis of Elliptic Boundary Value Problems and Eigenvalue Problems. PhD thesis, Massachusetts Institute of Technology, 1971.
[21] R. Courant, K. Friedrichs, H. Lewy. Über die partiellen Differenzengleichungen der mathematischen Physik. *Mathematische Annalen.* 1928; 100:32–74.
[22] R. Bhatia. *Matrix Analysis.* Graduate Texts in Mathematics 169; Springer-Science+Business Media, LLC; 1997.
[23] E. Fischer. Über quadratische Formen mit reellen Koeffizienten. *Monatshefte für Mathematik und Physik.* 1905; 16:234–249.
[24] R. Courant. Über die Eigenwerte bei den Differentialgleichungen der mathematischen Physik. *Mathematische Zeitschrift.* 1920; 7(1-4):1–57.
[25] H. Weyl. Das Asymptotische Verteilungsgesetz der Eigenwerte linearer partieller Differentialgleichungen (mit einer Anwendung auf die Theorie der Hohlraumstrahlung). *Mathematische Annalen.* 1912; 71:441–479.
[26] J. K. Merikoski, R. Kumar. Inequalities For Spreads Of Matrix Sums And Products. *Applied Mathematics E-Notes.* 2004; 4:150–159.
[27] C. R. Rao. Separation Theorems for Singular Values of Matrices and Their Applications in Multivariate Analysis. *Journal of Multivariate Analysis.* 1979; 9:362–377.
[28] R. A Horn, N. H. Rhee, W. So. Eigenvalue Inequalities and Equalities. *Linear Algebra and its Applications.* 1998; 44(1-3):29–44.
[29] L. Yu. Kolotilina. A Generalization of Weyl's Inequalities with Implications. *Journal of Mathematical Sciences.* 2000; 101(4):3255–3260.
[30] G. W. Stewart. *Matrix Algorithms.* Society for Industrial and Applied Mathematics; 1998.
[31] M. Argerami. Answer to "Singular value proofs", available at https://math.stackexchange.com/questions/244743/singular-value-proofs (last accessed 15/02/2018); 2012.
[32] J. A. Cottrell, T. J. R. Hughes, Y. Bazilevs. *Isogeometric Analysis: Toward Integration of CAD and FEA.* John Wiley & Sons, Ltd.; 2009.
[33] C. Adam, S. Bouabdallah, M. Zarroug, H. Maitournam. Stable time step estimates for NURBS-based explicit dynamics. *Computer Methods in Applied Mechanics and Engineering.* 2015; 295:581–605.
[34] L. Sirovich. Turbulence and the dynamics of coherent structures, part I-III.. *Quarterly of Applied Mathematics.* 1987; XLV(3):561–571.
[35] D. Amsallem, M. J. Zahr, C. Farhat. Nonlinear Model Order Reduction Based on Local Reduced-Order Bases. *International Journal for Numerical Methods in Engineering.* 2012; 92(10):891–916.
[36] K. Carlberg, M. Barone, H. Antil. Galerkin v. least-squares Petrov–Galerkin projection in nonlinear model reduction. *Journal of Computational Physics.* 2017; 330:693–734.
[37] P. R. Halmos. Properties of Spectra. In: A Hilbert Space Problem Book; Springer-Verlag New York, Inc.; 1982; 41–43.
[38] J. J. Sylvester. A Demonstration of the Theorem that Every Homogeneous Quadratic Polynomial is Reducible by Real Orthogonal Substitutions to the Form of a Sum of Positive and Negative Squares. *Philosophical Magazine.* 1852; 4(23):138–142.
[39] G. H. Golub, C. F. Van Loan. *Matrix Computations.* The John Hopkins University Press; Baltimore; 4th ed.2013.
[40] D. Serre. *Matrices - Theory and Applications.* Graduate Texts in Mathematics 216; Springer; 2nd ed. 2010.
[41] R. Everson, L. Sirovich. Karhunen-Loève procedure for gappy data. *Journal of the Optical Society of America A.* 1995; 12(8):1657–1664.
[42] Chaturantabut Saifon, Sorensen Danny C. Nonlinear Model Reduction via Discrete Empirical Interpolation. *SIAM Journal on Scientific Computing.* 2010; 32(5):2737–2764.
[43] K. Carlberg, R. Tuminaro, P. Boggs. Preserving Lagrangian structure in nonlinear model reduction with application

to structural dynamics. *SIAM Journal on Scientific Computing.* 2015; 37(2):153–184.

[44] F. Daim. Reduced order modelling for Crash Simulation. Doctoral workshop on model reduction in nonlinear mechanics; Mines ParisTech, Paris, France; 2017.

[45] K. Carlberg, C. Bou-Mosleh, C. Farhat. Efficient non-linear model reduction via a least-squares Petrov-Galerkin projection and compressive tensor approximations. *International Journal for Numerical Methods in Engineering.* 2011; 86(2):155–181.

[46] M. J. Zahr. Adaptive Model Reduction to Accelerate Optimization Problems Governed by Partial Differential Equations. PhD thesis, Stanford University, 2016.

[47] K. Carlberg, C. Farhat, J. Cortial, D. Amsallem. The GNAT method for nonlinear model reduction: Effective implementation and application to computational fluid dynamics and turbulent flows. *Journal of Computational Physics.* 2013; 242:623–647.

[48] P. Tiso, D. J. Rixen. Discrete Empirical Interpolation Method for Finite Element Structural Dynamics. In: Topics in Nonlinear Dynamics, pp. 203–212. Springer, 2013.

[49] S. Chaturantabut. Nonlinear Model Reduction via Discrete Empirical Interpolation. PhD thesis, Rice University, 2011.

[50] T. Chapman, P. Avery, P. Collins, C. Farhat. Accelerated mesh sampling for the hyper reduction of nonlinear computational models. *International Journal for Numerical Methods in Engineering.* 2017; 109(12):1623–1654.

[51] F. Ghavamian, P. Tiso, A. Simone. POD-DEIM model order reduction for strain softening viscoplasticity. *Computer Methods in Applied Mechanics and Engineering.* 2017; 317:458–479.

[52] A. Radermacher, S. Reese. POD-based model reduction with empirical interpolation applied to nonlinear elasticity. *International Journal for Numerical Methods in Engineering.* 2016; 107(6):477-–495.

[53] C. Farhat, T. Chapman, P. Avery. Structure-preserving, stability, and accuracy properties of the energy-conserving sampling and weighting method for the hyper reduction of nonlinear finite element dynamic models. *International Journal for Numerical Methods in Engineering.* 2015; 102(5):1077–1110.

[54] S. S. An, D. L. James. Optimizing Cubature for Efficient Integration of Subspace Deformations. *ACM Transactions on Graphics (TOG) - Proceedings of ACM SIGGRAPH Asia 2008.* 2008; 27(5).

[55] M. Géradin, D. J. Rixen. A "nodeless" dual superelement formulation for structural and multibody dynamics application to reduction of contact problems. *International Journal for Numerical Methods in Engineering.* 2016, 106:773—798.

[56] P. D. Lax, R. D. Richtmyer. Survey of the Stability of Linear Finite Difference Equations. *Communications on Pure and Applied Mathematics.* 1956; IX:267–293.